\documentclass[hidelinks]{article}
\usepackage[utf8]{inputenc}

\usepackage{graphicx,multirow,amsfonts,subfigure,hyperref,color}
\usepackage[normalem]{ulem}
\usepackage{amssymb,amsmath,amsfonts}
\usepackage{bm}
\usepackage{booktabs}

\usepackage{authblk}

\usepackage{algorithm}
\usepackage{algpseudocode}

\def\half{\frac{1}{2}}

\newcommand{\eq}[1]{\begin{equation}\label{#1}}
\newcommand{\en}{\end{equation}}
 
\def\Diag{\texttt{Diag}} 
\def\inv{^{-1}}%
\def\nref#1{(\ref{#1})}
\def\RR{\mathbb{R}}

\def\Span#1{\text{Span}{(#1)}}
\DeclareMathOperator{\myspan}{Span}

\DeclareMathOperator{\dist}{dist}

\DeclareMathOperator{\grad}{grad}
\DeclareMathOperator{\diag}{diag}

\usepackage{amsmath,amsfonts,amssymb,amsthm}

\newcommand{\bigO}{\mathcal{O}}
\newtheorem{theorem}{Theorem}

\newtheorem{corollary}[theorem]{Corollary}
\newtheorem{lemma}[theorem]{Lemma}
\newtheorem{proposition}[theorem]{Proposition}

\def\tla{\tilde \lambda}%
\def\tlu{\tilde u}%
\def\scX{{\cal X}} 

\DeclareMathOperator{\Tr}{Tr}
\DeclareMathOperator{\Gr}{Gr}

\def\betab{\vspace{-0.1in}\begin{tabbing}
xxx\=xxx\=xxxx\=xxx\=xxxxxxxxxxxxxxxxxxxx\=xxxx\=xxxx\= \kill} 
\def\entab{\end{tabbing}\vspace{0.1in}}

\newtheorem{algor}{\rlap{\rule[-5pt]{\textwidth}{1pt}}
\sc \color{red} {ALGORITHM :}}

\graphicspath{{./FIGS/}}

\title{Gradient-type subspace iteration methods for the
  symmetric eigenvalue problem}

\author[1]{Foivos Alimisis}
\author[2]{Yousef Saad}
\author[3]{Bart Vandereycken}
\affil[1,3]{Department of Mathematics, University of Geneva}
\affil[2]{Department of Computer Science and Engineering, University of Minnesota}

\date{} 

\begin{document} 

\maketitle 

\begin{abstract}
  This paper explores variants of the subspace iteration algorithm for computing
  approximate invariant subspaces. The standard subspace iteration approach is revisited   and new  variants that  exploit gradient-type
  techniques combined with a Grassmann manifold viewpoint are developed.  A 
  gradient method as well as a nonlinear conjugate gradient technique are described.
  Convergence of the gradient-based algorithm is analyzed and a few numerical
  experiments are reported, indicating that the proposed algorithms are sometimes superior to standard algorithms. This includes the Chebyshev-based subspace iteration and the locally optimal block conjugate gradient method, when compared in terms of number of matrix vector products and computational time, resp. The new methods, on the other hand, do not require
  estimating optimal parameters. An important contribution of this paper to achieve this good performance is the accurate and efficient implementation of an exact line search. In addition, new convergence proofs are presented for the non-accelerated gradient method that includes a locally exponential convergence if started in a  $\mathcal{O}(\sqrt{\delta})$ neighbourhood of the dominant subspace with spectral gap $\delta$.
  \end{abstract}

\textbf{AMS:} 15A69, 15A18


\section{Introduction}\label{sec:intro} 
When considering the many sources of large eigenvalue problems in numerical
linear algebra, one often observes that the actual underlying problem it to compute
an invariant subspace. In these cases, the 
eigenvalues and eigenvectors are often a by-product of the 
computations and they are not directly utilized.  For example, one of the most
common calculations in data science consists of performing a dimension
reduction which extracts a subspace that provides a good approximation of the original data 
in that not much information is lost when we project the original problem into this
low-dimensional space.  This  projection often results in better accuracy
since the information that is shed out corresponds to noise.
Another example is in  electronic structure calculations where 
an important class of algorithms called `linear scaling methods'
are entirely based on the eigenprojector on the subspace associated with the
`occupied states'. This projector is available through any orthonormal basis of the
invariant subspace and here again eigenvectors and eigenvalues are not explicitly
needed,  resulting in methods that scale linearly with the number of particles. 

While this distinction is often blurred in the literature,  
a number of articles did specifically deal with the problem of computing an
invariant subspace by expressing it  in terms of computing objects on the Grassmann manifold.
Thus, the well-known  article~\cite{eas:99}  addressed  general
optimization problems on manifolds, including eigenvalue problems as a special case.
A number of other papers also adopted, explicitly or implicitly, a matrix manifold viewpoint
for computing invariant subspaces, see, 
e.g., \cite{AbsilAl02,chatelin84,SamehTong00,BouchardAl18,alimisis2021distributed,ahn2021riemannian,alimisis2022geodesic}, among others. In many of these contributions, a Newton-type approach is advocated to solve the resulting equations.
Newton's method typically requires solving linear systems, which in this context are Sylvester equations, and this can be quite
expensive, or even impractical in some situations.

In view of the above observation, we may ask   whether or not a gradient approach can yield an effective alternative to standard implementations of subspace iteration.
Subspace Iteration (SI) displayed in Algorithm \ref{alg:SubsIt} with $p_k(t)=t$ computes the dominant  invariant subspace of a matrix and some of its known advantages are
highly desirable in a number of situations. SI is a block form of the power method and as such it  is
rather simple to implement. It is also known for its  robustness properties. For example, it is
resilient  to small  changes in the matrix during iteration,  an important attribute that is 
not shared by Krylov subspace methods. This particular feature is appealing in many practical instances as in, e.g, 
subspace tracking \cite{MoonenAl92,HoDeLaSuVan.04,doukopoulos2008fast,stewart1992updating,comon1990tracking},
or in electronic structure calculations \cite{YZhouAl14,chebfsi}.

This paper considers variations  of subspace iteration that are grounded in a gradient-type method on
the Grassmann manifold.
A gradient and a (nonlinear) conjugate gradient approach are described that both share
the  same advantages as those of classical SI. However, the proposed methods are based on a direct optimization of the partial trace of the matrix in the subspace. The convergence of the gradient algorithm will be studied theoretically.

We point out that other gradient-descent type methods that employ a Grassmannian
viewpoint have been recently developed in 
\cite{alimisis2021distributed,alimisis2022geodesic}.  These methods differ from those of this paper 
in that they aim at following a geodesic on the manifold by exploiting a
Riemannian structure. No such attempt is made in the methods proposed
herein. Instead, a simple gradient descent (or ascent) approach with (exact) line search
is applied and a simple normalization (retraction) step is added to ensure the orthonormality of the
basis of the new subspace.

The first part of the paper discusses classical versions of subspace iteration.
The second part develops line search techniques combined with gradient descent-type 
approaches.  A conjugate gradient approach is also presented. We end with numerical experiments.

\section{Background and notation}
This section begins with  providing  some background on
invariant subspaces and then defines the  notation to be used throughout the paper.

\subsection{Invariant subspaces}
Given an $n \times n $ matrix $A$, a subspace $\scX$ of $\mathbb{R}^n$ of dimension $p$ is called
invariant with respect to $A$ \emph{iff}:
\eq{eq:Inv}
A \scX \subseteq \scX .
\en
This can be expressed in matrix form by using a basis $X \in \RR^{n \times p}$ of
$\scX$. In this case $\scX$ is invariant \emph{iff} there exists a matrix
$\Lambda \in \RR^{p \times p}$ such that 
\eq{eq:Inv2}
A X = X \Lambda . 
\en
This second definition depends on a basis which is not  unique. Herein lies a conundrum
that is encountered in this context. We need a (non-unique) basis for
computations and expressing equations and equalities; however, the original definition
\nref{eq:Inv} does not require a basis.

A number of computational tasks deal specifically with invariant subspaces. 
The most common of these is probably just to \emph{compute} an invariant subspace as
represented by, e.g., an orthonormal basis. 
Sometimes, the task is to \emph{update} the subspace rather than compute it from
scratch. This happens  for example when solving the Kohn--Sham equation, see e.g.,
\cite{SaChSh-Nano} in electronic
structure calculations where at each iteration of the Self-Consistent Field (SCF) method the Hamiltonian changes
slightly and it is necessary to update an already computed approximate invariant subspace for
the previous Hamiltonian.
There are also numerous applications in signal processing, where the problem is to
\emph{track} an invariant  subspace of a sequence of matrices,
see e.g., \cite{stewart1992updating,doukopoulos2008fast,perry2010minimax,HoDeLaSuVan.04} 
a problem that is somewhat related to subspace updating problem.

Another problem that is often encountered is to (inexpensively)  estimate the dimension 
of some invariant subspace. Thus, the approximate rank or
numerical rank of some data matrix can be needed
in order to determine the proper dimension required for an adequate
`dimension reduction', or for subspace tracking~\cite{perry2010minimax}. 
This  numerical rank can be
determined as the dimension of the (near) invariant subspace corresponding to singular values
larger than a certain threshold $\epsilon$, see, e.g., \cite{UbaruYSrank16}. Another problem in signal processing, is  to
find a subspace that is simultaneously a near-invariant subspace for a set
of matrices. 
A common characteristic of  the examples just mentioned is that they all deal with
invariant subspaces -- but they do not require eigenvalues and vectors explicitly.

\subsection{Notation and assumptions}
For the remainder of the paper we will restrict our attention to the  case
when $A$ is real symmetric and positive definite. The positive-definiteness assumption is made for theoretical reasons and without loss of generality since the iterates produced by the algorithms in this paper are invariant to the transformation $A+cI$ for real $c$. In addition, we will consider
the problem of computing the invariant subspace associated with the largest $p$ eigenvalues, which we refer to as the $p$th ``dominant'' subspace. In case the subspace that corresponds to the smallest $p$ eigenvalues is sought, the algorithms can be applied to $-A$.

Given an $n \times n$  symmetric real matrix $A$, 
we denote by $\lambda_1 \geq \lambda_2 \geq \cdots \geq \lambda_n$ its 
eigenvalues counted with their multiplicities.

A common method for finding the  $p$th dominant subspace of $A$ consists
of  minimizing the function
\eq{eq:phi}
 \phi(X) = - \tfrac{1}{2} \Tr(X^TAX) , 
\en
over the set of $n \times p$ matrices with orthonormal columns, i.e.,
such that\footnote{In the complex
  Hermitian case, we would minimize $-\frac{1}{2} \Tr(X^H A X)$ over matrices that satisfy $X^H X = I$ where
  $X^H$ is the transpose conjugate of $X$.} 
$X^T X = I$; see, e.g., \cite{SamehTong00,Sameh-trace-min,absilOptBook08,eas:99}.
In general, non-zero spectral gap will be assumed,
\eq{eq:gap}
\delta = \lambda_p - \lambda_{p+1} > 0.
\en
This condition implies that there exists a unique dominant subspace associated with the $p$ largest eigenvalues of $A$. Thus, minimizing the objective $\phi$ has a unique solution.

We denote by  $\Diag(A) $  the diagonal matrix whose diagonal entries
are the same  as those of  $A$.  The notation 
is overloaded by defining $\Diag(\alpha_1, \ldots, \alpha_n)$ to be
the diagonal matrix with diagonal entries $\alpha_1, \ldots, \alpha_n$. This dual use of
$\Diag(.)$ causes no ambiguity and is consistent with common usage as, for example, in \textsc{Matlab}.

\subsection{Subspace iteration}\label{sec:SI}
Given some initial subspace with a basis $X_0 \in \RR^{n \times p}$,
the \emph{classical} subspace iteration
algorithm is nothing but  a Rayleigh--Ritz projection method onto
the subspace spanned by $X_k = A^k X_0$. That is,
we seek an approximate  eigenpair 
$\tla, \tlu$ where $\tla \ \in \ \RR$ and $\tlu \ \in \ \Span{X_k}$,
by requiring that $(A - \tla I ) \tlu \perp \Span{X_k}$.
If $Q = [q_1, \ldots, q_m]$ is an orthonormal 
basis of $X_k$, and 
we    express the  approximate eigenvector as 
$\tilde u = Q \tilde y$, then this leads to
$Q^T(A - \tla I)\tlu =0$ which means that $\tla, \tilde y$ is a solution of the
projected eigenvalue problem
\eq{eq:rr1}
Q^T A Q \tilde y = \tilde \lambda \tilde y. 
\en

  A common and more effective
alternative is to define $X_k$ to be of the form $X_k = p_k(A) X_0$, in which
$p_k$ is some optimally selected polynomial of degree $k$.  The
original form discussed above which was described by Bauer~\cite{Bauer}, corresponds to using the monomial $p_k(t) \equiv t^k$ or the shifted monomial 
$p_k(t) = (t-c)^k$. Rutishauser later developed more advanced versions in which $p_k(t)$ was selected to be a shifted and scaled Chebyshev polynomial \cite{Rutishauser-subs,RITZIT}.

\begin{algorithm}
  \caption{\mbox{\texttt{Subspace Iteration with polynomial filter}}$(A,X)$} \label{alg:SubsIt}
  \begin{algorithmic}[1]
    \State  \textbf{Start:} Select initial $X_0$ such that $X_0^T X_0 = I$ and polynomial $p_k$

    \For{$k = 0, 1, \ldots$} 
       \State  Compute $G := (I-X_k X_k^T)(A X_k)$.
    \If {$\| G \| < \textrm{tol}$}
    \State \Return
    \EndIf

    \State  Compute $\widehat X = p_k(A) X_k $
    \State Set $Q$ as the q-factor of the QR decomposition of $\widehat X$.
    \State Compute $C = Q^T A Q$.
    \State Diagonalize $C = U \Lambda_C U^T $ and set $X_{k+1} = Q U$ 
    \State Select a new polynomial $p_{k+1}$.    
\EndFor
\end{algorithmic}
\end{algorithm}

 \paragraph*{The optimal polynomial}\,
 An important issue when applying subspace iteration,
 is to select an optimal polynomial $p_k(t)$.  
 Assuming that we use a subspace of dimension $p$, the polynomial  is selected so as to enable the method to compute the eigenvalues $\lambda_1 \geq \lambda_2 \geq \cdots \geq \lambda_p$.  The standard approach  \cite{Rutishauser-subs} for subspace iteration when computing the dominant subspace is to use the polynomial
 \begin{equation}\label{eq:optimal_poly}
p_k(t) \equiv C_k((t-c)/h) \quad \text{where $c = (\lambda_{p+1}+\lambda_n)/2$ and $h = (\lambda_{p+1}-\lambda_n)/2$}.
 \end{equation}
Here, $C_k(t)$ is the Chebyshev polynomial of the first kind of degree $k$. Remark that $c$ is the middle of the interval $[\lambda_n, \ \lambda_{p+1}]$
 of unwanted eigenvalues, and $h$ is half the width of this interval.

 The polynomial above is found from an analysis of subspace iteration that reveals that, 
 when considering each eigenpair $\lambda_i, u_i$ for $i\le p$,
 the method acts as if the other eigenvalues $\lambda_j$ for $j \le p, j\ne i$ are not present.
 In other words, the polynomial is selected with the aim to minimize the maximum 
 $p_k(\lambda_\ell)/p_k(\lambda_i)$ for $\ell >p$ over polynomials of degree  $k$.
 This is  then relaxed to minimizing the maximum of
 $p_k(t)/p_k(\lambda_i)$ for $t \in \ [\lambda_n, \lambda_{p+1}]$ over
  polynomials of degree  $k$.
  The optimal polynomial is $C_k((t-c)/h)/C_k((\lambda_i-c)/h)$ where the denominator is
  for scaling purposes.
  In practice $\lambda_{p+1}, \ldots, \lambda_n$ are estimated and the scaling
  is performed for $i=1$, i.e., with $\lambda_1$  which is also estimated
  \cite{Rutishauser-subs,RITZIT}.
  Note that this polynomial is optimal for \emph{each eigenvalue
    $\lambda_j, j\le p$, individually}.
  However, it is not optimal when considering the subspace as a whole:
  a few experiments will reveal that convergence can be much faster when we replace $\lambda_{p+1} $ in the
  definition of polynomial by $\lambda_{p+k} $ for some $k>1$. 
  Note that there does not seem to be a theoretical study of this empirical observation.

 \paragraph{Comparison with Krylov subspace methods} \, 
  It is known that when computing a small number of eigenvalues and vectors at
 one end of the spectrum,
 Krylov subspace methods such as the Lanczos method and
 Arnoldi's method are generally faster than methods based on subspace iteration. 
 Standard Krylov methods require only one starting vector and this can be seen
 as an advantage in that little data needs to be generated to start the
 algorithm.  A known disadvantage of standard Krylov methods is that they cannot, in theory,
 compute more than  one eigenvector associated with a multiple eigenvalue.
 More importantly, for many applications it can be a disadvantage to start
 with one vector only.  Indeed, there are
 applications where a subspace must be computed repeatedly with a
 matrix that changes slightly from step to step. At the start of a new
 iteration we have available the whole (orthonormal) basis from the
 previous iteration which could be exploited. However, since Krylov
 methods start with only one vector this is not possible.

 On the other hand, the subspace iteration algorithm is perfectly suitable for the situation
 just  described: When a computation with a new matrix starts, we can take as
 initial subspace the latest subspace obtained from the previous
 matrix. 
 This is the exact scenario encountered in electronic structure
 calculations~\cite{YZhouAl14,chebfsi}, where a subspace is to be computed at
 each SCF iteration. The matrix changes at each SCF
 iteration and the changes depend only on the invariant
 subspace obtained at the previous SCF step.


  \section{Invariant subspaces and the Grassmannian perspective}\label{sec:grass} 
An alternative to using subspace iteration for computing a dominant
invariant subspace is to consider a method whose goal is to optimize the
objective function $\phi(X)$, defined in \eqref{eq:phi}, 
over all matrices $X$ that have orthonormal columns. 
This idea is of course not new. For example,  it is  the main ingredient exploited in the TraceMin
algorithm~\cite{SamehTong00,Sameh-trace-min}, a  method designed 
 for computing an invariant subspace associated with smallest 
eigenvalues for standard and generalized eigenvalue problems. 

A key observation in the definition~\nref{eq:phi} is that $\phi(X) $ is
invariant upon orthogonal transformations. In other words if $W$ is a $p \times p$
orthogonal matrix then, $\phi( X W) = \phi(X)$.  Noting that two orthonormal bases
$X_1$ and $X_2$ of the same subspace of dimension $p$ are related by $X_2 = X_1 W$ where $W$ is
some orthogonal $p \times p $ matrix, this means that the objective function
\nref{eq:phi} depends only on the subspace spanned by $X$ and not the particular
orthonormal basis $X$ employed to represent the subspace.  This in turn suggests
that it is possible, and possibly advantageous, to seek the optimum solution in
the Grassmann manifold~\cite{eas:99}.  Recall, from e.g., \cite{eas:99}, that
the Stiefel manifold is the set  
\eq{eq:Stnp} \textnormal{St}(n,p) = \{ X \ \in \ \RR^{n \times p} \colon \
X^T X = I \} , \en while the Grassmann manifold is the quotient manifold
\eq{eq:Grnp} \Gr(n,p) = \textnormal{St}(n,p) / \textnormal{O}(p) \en where $\textnormal{O}(p)$ is the orthogonal group of
unitary $p \times p$ matrices.  Each point on the manifold, one of the
equivalence classes in the above definition, can be viewed as a subspace of
dimension $p$ of $\RR^n$. In other words,
\begin{equation}\label{eq:grnk}
    \Gr(n,p)=\lbrace \mathcal{X} \subseteq \mathbb{R}^n \colon \mathcal{X} \hspace{1mm}  \text{is a subspace of} \dim(\mathcal{X})=p \rbrace.
\end{equation}
An element of $\Gr(n,p)$  can be indirectly represented by a basis
$X \in \textnormal{St}(n,p)$ modulo an orthogonal transformation and so we denote it by
$[X]$, keeping in mind that it does not matter which member $X$ of the
equivalence class is selected for this representation.

With this Riemannian viewpoint in mind, one can try to minimize $\phi(X)$ over the Grassmann manifold using one
of the many (generic) Riemannian optimization algorithms for smooth optimization. One has, for example, first-order and second-order algorithms that generalize
gradient descent and Newton's iteration, resp. We refer to the foundational article~\cite{eas:99} and the textbook~\cite{absilOptBook08} for detailed explanations. 
Despite its algorithmic simplicity, the Riemanian gradient method for this problem has not been treated in much detail beyond the basic formulation and its generically applicable theoretical properties. In the next sections, our aim is to work this out in detail and show in the numerical experiments that such a first-order method can be surprisingly effective and competitive.

From now on, our cost function is
\begin{equation}
    \phi([X])=-\tfrac{1}{2} \Tr(X^T A X), \quad \text{with} \ X^T X=I.
\end{equation}
For simplicity of notation though, we will be writing $\phi(X)$ instead of $\phi([X])$, using $X$ to represent the subspace spanned by it.

\subsection{Gradient method on Grassmann} \label{sec:gradls}

In a gradient approach we would like to produce an iterate $X_{k+1} \in \mathbb{R}^{n \times p}$ starting from $X_k \in \mathbb{R}^{n \times p}$ following a rule of the form
\eq{eq:tilY}
X_{k+1} = X_k - \mu \grad \phi(X_k),
\en
where the step size $\mu>0$ is to be determined by some line search. The 
direction opposite to the gradient is a direction of decrease for the objective
function $\phi$. 
However, it is unclear what value of the step $\mu$
yields the largest decrease in the value of $\phi$.
This means that some care has to be exercised in the search for the optimal
$\mu$. 
A gradient procedure may be appealing if a good approximate solution is already
known, in which case, the gradient algorithm may provide a less expensive
alternative to one step of the subspace iteration in Alg.~\ref{alg:SubsIt}.

For a Riemannian method defined on a manifold, the search direction (here, $-\grad \phi(X_k)$) always lies in the tangent space of the current point (here, $X_k$) of said manifold. This makes sense since directions orthogonal to the tangent space leave the objective function constant up to first order in the step if the iterates are restricted to lie on the manifold. For our problem, an element of the Grassmann manifold is represented by a member
$X \in \textnormal{St}(n,p)$ of its class. The tangent space of the Grassmann manifold at this
$[X]$ is the set of matrices $\Delta \in \RR^{n\times p}$ satisfying the
following orthogonality relation, see \cite{eas:99}:
\eq{eq:Tang} X^T \Delta = 0.
\en
The Riemannian gradient of \nref{eq:phi} at $[X]$ is 
\eq{eq:gradPhi}
\grad \phi(X) = -\Pi A X \equiv    -(A X -X C) , 
\en
with the orthogonal projector $\Pi = I-X X^T$, and the projected matrix $C = X^TAX$.

  Even though $\grad \phi(X_k)$ is in the tangent space
(and a direction of decrease for $\phi$), we are not interested in
$X_{k+1}$ per se but in the subspace that it spans. In particular,
since we use orthonormal bases to define  the value of $\phi$ on the manifold,
we will need  to `correct' the non-orthogonality of the update \nref{eq:tilY}
when considering $\phi$.  This will be discussed shortly.
For now we establish a few simple relations.

\begin{figure}[H]
  \centerline{\includegraphics[width=0.59\textwidth]{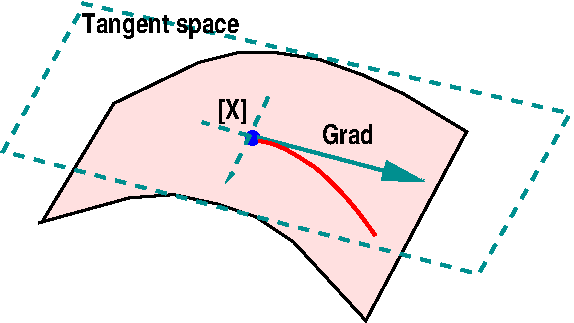} }
  \caption{Illustration of the line search and the tangent space}\label{fig:tangent}
\end{figure}

For simplicity we denote $X:=X_k$ an orthonormal basis of the current iterate, $\tilde X:=X_{k+1}$ a (probably non-orthonormal) basis of the new iterate and $G:=\grad \phi(X)$ the gradient direction. Then a step of the gradient method satisfies
$\tilde X = X - \mu G $ and we have
\eq{eq:phitil}
\phi(\tilde X) = \phi(X) - \mu \Tr  ((AX) ^T \Pi (AX)) -
\frac{\mu^2}{2} \Tr  ((AX) ^T \Pi A \Pi (AX)) .
\en
We also have the following relations
\begin{eqnarray}
  (AX) ^T \Pi (AX) &=&    -(AX) ^T  G = -(G^T (AX))^T = -G^T (AX)   \label{eq:GTAY}\\ 
                   &=&  (AX) ^T \Pi^T  \Pi (AX) = G^T G \label{eq:GTG}
\end{eqnarray}
where the second equality exploits the fact that  $\Pi $ is an orthogonal projector.

Thus, the coefficient of $\mu$ in the right-hand side of 
\nref{eq:phitil}  is nothing but $\|G \|_F^2$ and, therefore,
as expected,  the direction of $G$ is a descent direction: for small enough $\mu$,
$\tilde X$ will be close to orthonormal, and
regardless of the value of the trace in the last term, we would get a
decrease of the objective function $\phi$.
This will be the case unless we have already reached a critical point where $G=0$.

When looking at \nref{eq:phitil} it may appear at first that
when $A$ is SPD, it is possible to increase the value of $\mu$  arbitrarily and
decrease the objective function arbitrarily. This is clearly
incorrect because we have not yet adjusted the basis: we need to find the 
subspace spanned by $\tilde X$ and compute the related value of the objective function. 
In the following we address this issue by actually optimizing
the objective function on the manifold.

Observe that since  $X^T G = 0$ we have: 
\[
  \tilde X^T \tilde X = (X - \mu G)^T (X - \mu G)  = I + \mu^2 G^T G . 
\]
Let the spectral decomposition of $G^T G$ be
\eq{eq:Gdec} G^T G = V D_{\beta} V^T 
\en
and denote $\beta=\Diag(D_\beta)$ the eigenvalues. 
 We now define the diagonal matrix
\eq{eq:Psi3}
D_\mu \equiv (I + \mu^2 D_{\beta})^{1/2} .
\en

In order to make 
$\tilde X$ orthogonal 
without  changing its linear span,
we multiply it to the right by $V D_\mu \inv V^T$. This way, we obtain the matrix
\eq{eq:Ynew}
X(\mu) = \tilde X V D_\mu\inv V^T = (X - \mu G) V D_\mu\inv V^T. 
\en
that depends on the step $\mu$ and is easily seen to be orthonormal,
\begin{eqnarray*}
  X(\mu)^T  X(\mu) = D_\mu\inv V^T (I + \mu^2  G^T G ) V  D_\mu\inv = I. 
\end{eqnarray*}

While it is tempting to remove the $V^T$ in~\eqref{eq:Ynew} as this does not change the linear span, it is useful to keep it. The normalization is only then equivalent to the polar factor of $X - \mu G$. In the context of optimization on manifolds, this so-called retraction has many nice properties. In particular, $X(\mu)$ is a best approximation of $X - \mu G$ in the set of orthonormal matrices. In addition, this retraction has an easy vector transport that is invariant to the choice of representative in the subspace, which will be important later in~\S\ref{sec:CG} when we discuss the acceleration of the gradient method.

\subsection{Efficient line search}\label{sec:opt}
We can  now tackle the issue of determining  the optimal $\mu$. If we set
\eq{eq:YuGu}
X_v = X V, \qquad G_v = GV ,
\en
then from \nref{eq:GTAY}--\nref{eq:GTG}  we get the relation
$ G_v^T A X_v = -G_v ^T G_v$. In addition,
  note that $G_v ^T G_v = V^T G^T G V = D_{\beta}$.
With these relations we can now show:
  \begin{eqnarray}
  \phi (X(\mu))
  &=& -\tfrac{1}{2} \Tr  (VD_\mu\inv V^T (X-\mu G)^T A (X-\mu G) V D_\mu\inv V^T) \nonumber \\
  &=& -\tfrac{1}{2}  \Tr  (D_\mu\inv (X_v-\mu G_v)^T A (X_v-\mu G_v) D_\mu\inv ) \nonumber \\
  &=& -\tfrac{1}{2} \Tr  \left( D_\mu^{-2} \left( X_v^T A X_v + 2 \mu (G_v^T G_v )  
     + \mu^2 ( G_v^T A G_v )\right) \right)  \nonumber \\
  &=& -\tfrac{1}{2} \Tr   \left(   \left(I+\mu^2 D_{\beta} \right)^{-1} 
      \left( X_v^T A X_v + 2 \mu \ D_\beta 
        + \mu^2  \ G_v^T A G_v \right) \right) .     \label{eq:ratmu} 
  \end{eqnarray}
  We will simplify notation by introducing the  diagonal matrices:
  \begin{align}
    D_{\alpha}  &= \Diag(\alpha_1, \ldots, \alpha_p)  \quad  \text{with} \quad\alpha_i = ( X_v^T A X_v)_{ii}, \label{eq:diags}  \\ 
    D_\gamma &= \Diag(\gamma_1, \ldots, \gamma_p) \quad  \text{with} \quad \gamma_i = ( G_v^T A G_v)_{ii} . 
       \label{eq:diags2} 
  \end{align}

  If we call $u_i$ the left singular vector of $G$ associated with
$\sqrt{\beta_i}$ then  we get the useful relation
\eq{eq:gamRel}  \gamma_i \equiv v_i^T G^T A G v_i = \beta_i u_i^T A u_i . \en

  Observe that when $D$ is a diagonal matrix and $C$ is arbitrary, then
  $\Diag(D C) = D \ \Diag(C)$. Therefore,   \nref{eq:ratmu} simplifies to:
  \eq{eq:ratmu1}
  \phi (X(\mu))  = - 
  \tfrac{1}{2} \Tr   \left(  \left(I+\mu^2 D_{\beta} \right)^{-1} 
  \left( D_\alpha + 2 \mu D_{\beta}  + \mu^2 D_{\gamma} \right) \right) \ .
  \en   
  This is a rational function that is the sum of $p$ terms corresponding to the $p$
  diagonal entries of the matrix involved in \nref{eq:ratmu1}:
  \eq{eq:phinew} 
    \phi (X(\mu)) = - \half \sum_{i=1}^p  
    \frac{ \alpha_{i}\ + 2 \beta_i  \mu +  \gamma_{i} \mu^2 }
    {1 +  \beta_i \mu^2} . 
  \en

   When $\mu \to \infty$ each term $\frac{ \alpha_{i}\ + 2 \beta_i  \mu +  \gamma_{i} \mu^2 }
    {1 +  \beta_i \mu^2}$ will \emph{decrease}
    to its limit $\gamma_i /\beta_i$.
  The derivative of $\phi(X(\mu))$ satisfies
  \eq{eq:dqnew}
    \frac{ d \phi(X(\mu))}{d \mu} = -\sum_{i=1}^p \
    \frac{ \beta_i\ + (\gamma_{i}\ -  \alpha_{i}\ \beta_i) \mu -      \beta_i^2 \mu^2 }
    {(1+ \beta_i \mu^2)^2} \ .
    \en
    This derivative is the negative sum of $p$  branches each associated with a
    diagonal entry of the matrix of which the trace is taken in the above
    equation.      The numerator
    $\beta_i + (\gamma_{i} -  \alpha_{i} \beta_i) \mu -      \beta_i^2 \mu^2  $
     of each branch has the shape of an inverted parabola and has a negative and a positive root.
     Therefore, the derivative \nref{eq:dqnew} is nonpositive at zero\footnote{It is equal to
       $-\sum \beta_i = -\|G\|_F^2$} and
     as $\mu$ increases away from the
     origin, each of the branches will have a negative derivative.
    The derivative remains negative until $\mu$ reaches the second root which is
    \eq{eq:xii}
      \xi_{i}  =   \frac{   (\gamma_{i}\ -  \alpha_{i} \beta_i) + 
        \sqrt{ (\gamma_{i}\ -  \alpha_{i} \beta_i)^2 +
          4 \beta_i^3 } }
      { 2 \beta_i^2 }  \ > \  0  . 
       \en  
       
  Let $\xi_{min} = \min_i \{ \xi_i \}$ and $\xi_{max} = \max_i \{  \xi_i \}$.
  Clearly all branches of \nref{eq:phinew}, and therefore also their sum,
  will decrease in value
  when $\mu $ goes from zero to $\xi_{min}$. Thus, the value of the objective function \nref{eq:phinew} will decrease. Similarly,
  when $\mu $ increases from $\xi_{max}$  to infinity, the objective function
  \nref{eq:phinew} will increase. The minimal value of \nref{eq:phinew}  with respect to
  $\mu$ can therefore be determined by seeking the minimum in the interval  $[\xi_{min}, \ \xi_{max}]$. Since both $\phi$ and its derivative are available, this can be done efficiently by any standard root finding algorithm.

  The algorithm to get the optimal value for $\mu$ is described in
  Algorithm \ref{alg:getmu}. To obtain accurate solutions, some care is required in the numerical implementation due to floating point arithmetic. We explain this  in more detail in Section~\ref{sec:numerical_implementation_details}.

\begin{algorithm}
  \caption{$\mbox{\texttt{Riemannian Gradient Descent}}(A,X)$}\label{alg:gradA1} 
  \begin{algorithmic}[1]
    \State  \textbf{Start:} Select initial $X_0$ such that $X_0^T X_0 = I$.

    \For{$k = 0, 1, \ldots$} 
    \State  Compute $G := \grad \phi(X_k) = -(A X_k - X_k C_k)$ with $C_k = X^T_k A X_k$.
    \If {$\| G \| < \textrm{tol}$}
    \State \Return
    \EndIf

    \State Diagonalize  $G^T G = V D_{\beta} V^T$.

    \State Compute $D_{\alpha}, D_{\gamma}$ from \nref{eq:diags} with $X=X_k$.

    \State Compute $\mu$ as the (approximate) minimizer \nref{eq:phinew} using  \texttt{Get\_Mu}.

    \State Compute $X_{k+1}$ as the polar factor of $X_k - \mu G$ like in~\eqref{eq:Ynew}. 
\EndFor
\end{algorithmic}
\end{algorithm}

\begin{algorithm}
  \caption{$\mu_{out} = \texttt{Get\_Mu} (D_{\alpha}, D_{\beta},
    D_{\gamma})$}\label{alg:getmu} 
  \begin{algorithmic}[1]
    \State  \textbf{Input:} Diagonal matrices 
    $D_{\alpha}, D_{\beta}, D_{\gamma}$ of \nref{eq:diags}.

    \State  Compute smallest root
    $\xi_{min} $  and largest root $\xi_{max} $ among the roots $\xi_i$
    of \nref{eq:xii}

    \State Compute an approximation $\mu_{out}$ of the minimum of $\phi$ on 
    $[ \xi_{min}, \ \xi_{max}] $ by safe-guarded root finding on~\eqref{eq:phinew}.

    \State \textbf{Return:} value $\mu_{out}$    
\end{algorithmic}
\end{algorithm}

\section{Convergence of the gradient method}\label{sec:conv}
We prove that the gradient method from Algorithm~\ref{alg:gradA1} converges globally to a critical point, that is, where the Riemannian gradient is zero. This result is valid for any initial iterate $X_0$ but it does not give a linear rate of convergence. When $X_0$ is close to the dominant subspace, we also prove a linear (exponential) rate  of convergence of the objective function. The closeness condition depends on the spectral gap $\delta$ of the dominant subspace but only as $O(\sqrt{\delta})$. This result seems to be new.

\subsection{Global convergence of the gradient vector field}

We 
examine the expression \eqref{eq:phinew} in order to obtain a useful lower bound. We first rewrite \eqref{eq:phinew} as follows:
  \begin{align} 
  \phi (X(\mu))
  &= -\half \sum_{i=1}^m  
    \frac{ \alpha_{i} (1+  \beta_i \mu^2) - \alpha_i  \beta_i \mu^2 
      + 2 \beta_i  \mu +  \gamma_{i} \mu^2 }
    {1 +  \beta_i \mu^2} \nonumber \\
  & = -\half \sum_{i=1}^m  \alpha_i -  \half \sum_{i=1}^m 
    \frac{ 2 \beta_i  \mu +  (\gamma_{i} - \alpha_i \beta_i)  \mu^2 }
    {1 +  \beta_i \mu^2}  \ . \label{eq:phiconv0} 
  \end{align}
  The first sum on the right-hand side
  is just the objective function before the update, that is, 
  the value of $\phi$ at the current iterate $X(0)=X$. The second sum depends on the step $\mu$ and thus represents what
  may be termed the `loss' of the objective function for a given $\mu$.

  \begin{lemma}\label{lem:conv1}
    Define $L \equiv \lambda_{max}(A) - \lambda_{min}(A)$.
    Then for any given $\mu \ge 0$
    the `loss' term (2nd term in right-hand side of \eqref{eq:phiconv0})
    satisfies
    \eq{eq:gainLem}
    \half \sum_{i=1}^m
    \frac{ 2 \beta_i  \mu +  (\gamma_{i} - \alpha_i \beta_i)  \mu^2 }
    {1 +  \beta_i \mu^2}  \ge 
  \frac{ (2   - L \mu) \mu  } 
  {2(1 +  \beta_{max}  \mu^2)}  \cdot  \| G \|_F^2 ,
  \en
  where $G=\grad \phi(X(0))$ and $\beta_{max}  = \max \beta_i$.
  \end{lemma}
\begin{proof}   
  We exploit \eqref{eq:gamRel} and set $\tau_i = u_i^T A u_i$ in order
  to  rewrite the term
  $  \gamma_{i}\ -  \alpha_{i}\ \beta_i $ in the numerator as
  $  \gamma_{i}\ -  \alpha_{i}\ \beta_i = (\tau_i - \alpha_i) \beta_i$.   From~\eqref{eq:YuGu} and~\eqref{eq:diags2}, we have $\alpha_i = x_i^T A x_i$ with $x_i = Xv_i$.  Hence, the
  term $\tau_i - \alpha_i \equiv u_i^T A u_i - x_i^T A x_i $ represents
  the difference between two Rayleigh quotients with respect to $A$ and
  therefore,    $\tau_i - \alpha_i \ge - L $. Thus the `loss` term
  satisfies 
    \eq{eq:gainPrf1}
  \half \sum_{i=1}^m  \frac{ 2 \beta_i  \mu +  (\gamma_{i} - \alpha_i \beta_i)  \mu^2 }
    {1 +  \beta_i \mu^2} \ge
  \half \sum_{i=1}^m  \frac{ 2   - L \mu } 
    {1 +  \beta_i \mu^2} \beta_i  \mu . 
    \en
    The denominators $1 +  \beta_i \mu^2$ can be bounded from above by
    $1 +  \beta_{max} \mu^2 $ and this will result in:
    \eq{eq:gainPrf2}
  \half \sum_{i=1}^m  \frac{ 2 \beta_i  \mu +  (\gamma_{i} - \alpha_i \beta_i)  \mu^2 }
    {1 +  \beta_i \mu^2} \ge
  \half \sum_{i=1}^m  \frac{ 2   - L \mu } 
    {1 +  \beta_{max}  \mu^2} \beta_i  \mu = 
  \frac{ (2   - L \mu) \mu  } 
    {2(1 +  \beta_{max}  \mu^2)}  \sum_{i=1}^m  \beta_i  .
    \en
    The proof ends by noticing that $\sum \beta_i = \| G \|_F^2$ due to~\eqref{eq:Gdec}.
\end{proof}

We now state a useful global upper bound for the biggest singular value of the Riemannian gradient. This result is proved in \cite[Lemma~4]{alimisis2022geodesic} using properties from linear algebra.\footnote{Note that in \cite{alimisis2022geodesic} the Rayleigh quotient is scaled by a factor $2$. This does not change the result, as both $G$ and $L$ are scaled accordingly.}
\begin{lemma}\label{lem:upper_bound_grad}
The spectral norm of the Riemannian gradient $G$ of $\phi$ satisfies
$ \| G \|_2  \leq L / 2$ at any point.
\end{lemma}

\begin{lemma}\label{lem:conv2} 
  If $\mu_{opt} $ is the optimal $\mu$ obtained from a line search at a given
  $X$,  then
\eq{eq:conv2}
      \phi (X(\mu_{opt})) \leq   
      -\half \sum_{i=1}^m \alpha_i -
    \frac{2}{5}  \frac{ \| G \|_F^2 } 
 { L}  \ . 
\en
  \end{lemma} 
  \begin{proof}
    The right-hand side \nref{eq:gainLem} is nearly minimized for
    $\mu_s =  1 / L$ so we consider this special value of $\mu$. We have
    \[
      \phi (X(\mu_{opt})) \le       \phi (X(\mu_{s}))  \le 
 -\half \sum_{i=1}^m \alpha_i   -    \frac{ (2   - L \mu_s) \mu_s  } 
 {2(1 +  \beta_{max}  \mu_s^2)}  \cdot  \| G \|_F^2 .
\]
The second inequality in the above equation follows from \eqref{eq:phiconv0}
and the previous Lemma \ref{lem:conv1}.
Calculating the right-hand side for $\mu_s = 1/L$ yields: 
\begin{equation*}
   \phi (X(\mu_{opt})) \leq -\half \sum_{i=1}^m \alpha_i - \frac{ \| G \|_F^2 }{2( L +  \beta_{max} /L)}.
\end{equation*}
Be Lemma \ref{lem:upper_bound_grad}, we have $\beta_{max} \leq \frac{L^2}{4}$ since $\beta_{max}$ is the biggest eigenvalue of $G^T G$. Plugging this into the last inequality we get the desired result.
\end{proof}

The property~\eqref{eq:conv2} in Lemma~\ref{lem:conv2} is known as a sufficient decrease condition of the line search. We can now follow standard arguments from optimization theory to conclude that (Riemannian) steepest descent  for the smooth objective function $\phi$ converges in gradient norm. 

  \begin{theorem}\label{prop:conv1}
    The sequence of gradient matrices $\grad \phi(X_k)$ generated by Riemanian gradient descent with exact line search
     converges (unconditionally) to zero starting from any $X_0$.
  \end{theorem}
  \begin{proof}
We will proceed by  avoiding the use of indices.  First, we observe that the
    traces of the iterates, that is, the consecutive values of $\phi(X(\mu_{opt}))$
    converge since they constitute a bounded decreasing sequence. Recall that
    the first term, that is, minus the half sum of the $\alpha_i$'s in the right-hand
    side of \eqref{eq:conv2}, is the value of the objective function at the
    previous iterate. Thus, the second term in \nref{eq:conv2} is bounded from
    above by the difference between two consecutive traces:
\begin{equation}\label{eq:diff fval step}
 0  \le   \frac{2}{5}  \frac{ \| G \|_F^2}{L}  \le 
  -\phi (X(\mu_{opt})) - 
  \half \sum_{i=1}^p \alpha_i = -\phi (X(\mu_{opt})) + \phi(X),
\end{equation}
and therefore it converges to zero. This implies that the sequence of gradients also converges to $0$.
\end{proof}

The bound of Lemma \ref{lem:conv2} can be used to prove some particular rate of convergence for the gradient vector field. This argument is again classical for smooth optimization. It is a slow (algebraic) rate but it holds for any initial guess. 
\begin{proposition}
    The iterates $X_k$ of Algorithm \ref{alg:gradA1} satisfy 
    \begin{equation*}
        \min_{k=0,...,K-1} \| \grad \phi(X_k)\|_F \leq \sqrt{\frac{5}{2}L (\phi(X_0)-\phi^*)} \frac{1}{\sqrt{K}},
    \end{equation*}
where $\phi^*$ is the minimum of $\phi$.
\end{proposition}

\begin{proof}
Since  $\phi^*$ is the minimum of $\phi$, it holds
    \begin{align}\label{eq:tele_sum}
        \phi(X_0)-\phi^* & \geq \phi(X_0)-\phi(X_K) = \sum_{k=0}^{K-1} (\phi(X_k)-\phi(X_{k+1})). 
    \end{align}
After some rearrangement, Lemma~\ref{lem:conv2} provides the bound
\begin{equation*}
   -\half \sum_{i=1}^m \alpha_i  - \phi(X(\mu_{opt})) = \phi(X_k)-\phi(X_{k+1}) \geq \frac{2}{5L} \|\textnormal{grad}\phi(X_k) \|_F^2.
\end{equation*}
Taking the sum of this inequality for $k=0,...,K-1$, we obtain the lower bound
\begin{equation*}
    \sum_{k=0}^{K-1} (\phi(X_k)-\phi(X_{k+1})) \geq K \frac{2}{5L} \min_{k=0,...,K-1} \| \grad \phi(X_k)\|_F^2.
\end{equation*}
Combining with~\eqref{eq:tele_sum} gives the desired result.
\end{proof}


\subsection{Local linear convergence}

The previous proposition establishes a global but slow convergence to a critical point. We now turn to the question of proving a fast (exponential or linear) rate to the dominant $p$-dimensional subspace $\mathcal{V}_\alpha=\textnormal{span}(V_{\alpha})$ of $A$. The result will only hold locally, however, for an initial guess $X_0$ sufficiently close to $\mathcal{V}_\alpha$. We therefore also assume a non-zero spectral gap $\delta=\lambda_{p} - \lambda_{p+1} >0$. We denote the global optimal value as $\phi^* = \phi(V_\alpha)$.

Let $T_{\mathcal{X}} \Gr(n,p)$ denote the tangent space of the Grassmann manifold $\Gr(n,p)$ at $\mathcal{X} \in \Gr(n,p)$ (represented by an orthonormal matrix). We  take the inner product between two tangent vectors in  $T_{\mathcal{X}} \Gr(n,p)$ as
$$
 \langle \Delta_1, \Delta_2 \rangle_{\mathcal{X}} = \Tr(\Delta^T_1 \Delta_2) \ \text{\ with\  $\Delta_1,\Delta_2 \in T_{\mathcal{X}} \Gr(n,p)$. }
$$
Here, $\Delta_1$ and $\Delta_2$ are tangent vectors of the same orthonormal representative $X$. Observe that the inner product is invariant to the choice of this representative
since the inner product of $\Bar{\Delta}_1=\Delta_1 R$ and $\Bar{\Delta}_2 = \Delta_2 R$
with orthogonal $R$, is the same as
$ \langle \Delta_1, \Delta_2 \rangle_{\mathcal{X}}$.
The norm induced by this inner product in any tangent space is the Frobenius norm, which is therefore compatible with our other theoretical results.

The Riemannian structure of the Grassmann manifold can be conveniently described
by the notion of \emph{principal angles} between subspaces. Given two
subspaces $\mathcal{X},\mathcal{Y} \in \Gr(n,p)$ spanned by the orthonormal matrices $X,Y$ respectively, the principal angles between
them are $0 \leq \theta_1 \leq \cdots \leq \theta_p \leq \pi/2$ obtained from
the SVD
\begin{equation}\label{eq:SVD_for_principal_angles}
    Y^T X=U_1 \cos \theta \ V_1^T
\end{equation}
where $U_1 \in \mathbb{R}^{p \times p}, V_1 \in \mathbb{R}^{p \times p}$ are orthogonal and $\cos \theta= \diag(\cos \theta_1,...,\cos \theta_p)$.

We can express the intrinsic distance induced by the Riemannian inner product discussed above as
\begin{equation}\label{eq:distance_with_Log_and_angles}
    \dist(\mathcal{X},\mathcal{Y)}=\sqrt{\theta_1^2+...+\theta_p^2}=\| \theta \|_2,
\end{equation}
where $\theta=(\theta_1, \ldots ,\theta_p)^T$.

The convexity structure of the Rayleigh quotient $\phi$ on the Grassmann manifold, with respect to the aforementioned Riemannian structure, is studied in detail in \cite{alimisis2022geodesic}. In the next proposition, we summarize all the important properties that we use for deriving a linear convergence rate for Algorithm \ref{alg:gradA1}.
In the rest, we denote subspaces of the Grassmann manifold by orthonormal matrices that represent them.

\begin{proposition}\label{prop:props of f}
\label{prop:quadratic_growth}
Let $0 \leq \theta_1 \leq \cdots \leq \theta_p < \pi/2$ be the principal angles between the subspaces $\mathcal{X}=\textnormal{span}(X)$ and $\mathcal{V}_\alpha=\textnormal{span}(V_{\alpha})$. The function $\phi$ satisfies
\begin{enumerate}
    \item $\phi(X)-\phi^* \geq c_Q \, \delta \, \textnormal{dist}^2(X,V_{\alpha})$ \qquad (quadratic growth)
    \item $\| \grad \phi(X) \|_F^2 \geq 4 \, c_Q \, \delta \, a^2(X) (\phi(X)-\phi^*)$ \qquad (gradient dominance)
    \item The eigenvalues of the Riemannian Hessian of $\phi$ are upper bounded by $ L$. This implies $\phi(X)-\phi^* \leq \tfrac{1}{2} L \textnormal{dist}^2(X, V_{\alpha})$  \qquad (smoothness)
    \item $\| \grad \phi(X) \|_2 \leq \tfrac{1}{2} L $ \qquad (cfr.~Lemma~\ref{lem:upper_bound_grad})
\end{enumerate}
where  $c_Q = 2/\pi^2$, $\delta = \lambda_p - \lambda_{p+1}$, $L = \lambda_{max}(A)-\lambda_{min}(A)$, and $a(X) = \theta_p / \tan \theta_p$.
\end{proposition}

The inequality mentioned in point 3 of Proposition \ref{prop:props of f} is derived by setting $Y=V_{\alpha}$ in Inequality (10) in \cite{alimisis2022geodesic}. Inequality (10) in \cite{alimisis2022geodesic} is directly implied by smoothness.

Next, we use these properties to prove an exponential convergence rate for the function values of $\phi$. In order to guarantee a uniform lower bound for $a(X_k)$ at the iterates $X_k$ of Algorithm \ref{alg:gradA1}, we need to start from a distance at most $\bigO(\sqrt{\delta})$ from the optimum.

\begin{proposition}
\label{thm:conv_function_val}
An iterate $X_{k+1}$ of Algorithm \ref{alg:gradA1} starting from a point $X_k$  satisfies
\begin{equation*}
    \phi(X_{k+1})-\phi^* \leq \left(1-\frac{8}{5} c_Q a^2(X_k) \frac{\delta}{L} \right) (\phi(X_k)-\phi^*).
\end{equation*}
\end{proposition}
\begin{proof}
The result follows simply by combining the bounds of Lemma \ref{lem:conv2} and Proposition \ref{prop:props of f} (2).
By Lemma \ref{lem:conv2}, we have
  \begin{equation*}
      \phi(X_{k+1})-\phi^* \leq \phi(X_k)-\phi^* - \frac{2}{5 L} \| \grad \phi (X_k)\|^2.
  \end{equation*}
  By the gradient dominance of $\phi$ in Proposition \ref{prop:props of f}, we have
  \begin{align*}
      \phi(X_{k+1})-\phi^* & \leq \phi(X_k)-\phi^* - \frac{8}{5} c_Q a^2(X_k) \frac{\delta}{L} (\phi(X_k)-\phi^*) \\  & \leq  
      \left( 1 -  \frac{8}{5} c_Q a^2(X_k) \frac{\delta}{L} \right) (\phi(X_k)-\phi^*).
  \end{align*}
  This provides the desired result.
  \end{proof}

The convergence factor in the previous theorem still involves a quantity $a(X_k)$ that depends on the iterate $X_k$ at step $k$. To get a convergence factor for all $k$ that only depends on the initial step, we need to bound $a(X_k)$ globally from below and independently of $k$. To that end, we need to restrict the initial guess $X_0$ to a radius $\mathcal{O}(\sqrt{\delta})$ away from the optimum. The reason for that is that, using Proposition \ref{thm:conv_function_val}, we can only show that function values do not increase. In order to obtain a bound for the distances of the iterates to the optimum (and thus also for $a(X_k)$), we need to use the quadratic growth condition of Proposition \ref{prop:props of f}. This leads to a loss of a  factor $\delta$ in the upper bound for the squared distances of the iterates to the optimum.

\begin{theorem}\label{cor:global_rate}
Algorithm \ref{alg:gradA1}, where $X_0$ is such that
$$\textnormal{dist}(X_0, V_{\alpha}) \leq \sqrt{\frac{2 c_Q \delta}{L}}, $$
produces iterates $X_k$ that satisfy 
\begin{equation*}
    \phi(X_k)-\phi^* \leq \left(1-c_Q\frac{2 \delta}{5 L} \right)^k (\phi(X_0)-\phi^*)
\end{equation*}
for all $k \geq 0$.
\end{theorem}

\begin{proof}
Recall that $a(X_k) = \theta_p / \tan \theta_p$ with $\theta_p$ the largest principal angle between $X_k$ and $V_{\alpha}$. 
By the result of Proposition \ref{thm:conv_function_val}, we have
 \begin{equation*}
    \phi(X_{k+1})-\phi^* \leq \left(1- \frac{8}{5} c_Q a^2(X_k) \frac{\delta}{L}  \right) (\phi(X_k)-\phi^*) \leq \phi(X_k)-\phi^*,
\end{equation*}
since $1- \frac{8}{5} c_Q a^2(X_k) \frac{\delta}{L} \leq 1$.
By induction, we can conclude that
\begin{equation*}
    \phi(X_k)-\phi^* \leq \phi(X_0)-\phi^*,
\end{equation*}
for all $k \geq 0$.

Then by quadratic growth and smoothness of $f$ in Proposition \ref{prop:props of f}, we have 
\begin{align*}
     \textnormal{dist}^2(X_k, V_{\alpha}) &\leq \frac{1}{c_Q \delta} (\phi(X_k)-\phi^*) 
     \leq \frac{1}{c_Q \delta} (\phi(X_0)-\phi^*) \\
     &\leq \frac{L}{2 c_Q \delta} \textnormal{dist}^2(X_0, V_{\alpha}) \leq 1,
\end{align*}
for all $k \geq 0$,
by the assumption on the initial distance between $X_0$ and $V_{\alpha}$.

By elementary properties of $\cos(x)$ and $x/\tan(x)$ and using~\eqref{eq:distance_with_Log_and_angles}, we have
\begin{equation*}
    a(X_k) \geq \cos(\theta_{p}(X_k,V_{\alpha})) \geq \cos(\textnormal{dist}(X_k,V_{\alpha})) \geq \cos(1) \geq \frac{1}{2}.
\end{equation*}
Plugging this in the result of Proposition \ref{thm:conv_function_val} and by an induction argument, we get the desired result.
\end{proof}

Finally, we present an iteration complexity for computing an approximation of the leading eigenspace via Algorithm \ref{alg:gradA1}. The $\Tilde{\bigO}$ notation hides non-leading logarithmic factors. This results is standard when a non-asymptotic convergence rate (like the one of Theorem \ref{cor:global_rate}) is available.
\begin{corollary}
Algorithm \ref{alg:gradA1} where $X_0$ satisfies the assumption of Thm.~\ref{cor:global_rate} computes an estimate $X_T$ of $V_{\alpha}$ such that $\textnormal{dist}(X_T,V_{\alpha}) \leq \epsilon$ in at most
\begin{equation*}
    T= \frac{5 \pi^ 2 L}{8 \delta} \log \frac{\phi(X_0)-\phi^*}{c_Q \varepsilon \delta} + 1 = \Tilde{\bigO}\left(\frac{L}{\delta } \log \frac{\phi(X_0)-\phi^*}{ \varepsilon } \right).
\end{equation*}
many iterates.
\end{corollary}

\begin{proof}
For $\textnormal{dist}(X_T,V_{\alpha}) \leq \epsilon$, it suffices to have
\begin{equation*}
    \phi(X_T)-\phi^* \leq c_Q \epsilon^2 \delta
\end{equation*}
by quadratic growth of $f$ in Proposition \ref{prop:props of f}. Using $(1-c)^k \leq \exp(-c k)$ for all $k \geq 0$ add $0 \leq c \leq 1$, Theorem \ref{cor:global_rate} gives that it suffices to choose $T$ as the smallest integer such that 
\begin{equation*}
    \phi(X_T)-\phi^* \leq \exp\left(- c_Q\frac{2 \delta}{5L} T \right) (\phi(X_0)-\phi^*) \leq c_Q \epsilon^2 \delta.
\end{equation*}
Solving for $T$ and substituting $c_Q = 4/\pi^2$, we get the required statement.
\end{proof}


\section{Accelerated gradient method}\label{sec:CG}

It is natural to consider an accelerated gradient algorithm as an
improvement to the standard gradient method. For convex quadratic functions on $\RR^n$, the best example is the conjugate gradient algorithm since it speeds up convergence significantly at virtually the same cost per step as the gradient method. In our case, the objective function is defined on $\Gr(n,p)$ and is no longer quadratic. Hence, other ideas are needed to accelerate. While there exist a few ways to accelerate the gradient method, they all introduce some kind of momentum term and compute a new search direction $P$ recursively based on the previous iteration. 

\subsection{Polak--Ribiere nonlinear conjugate gradients} 

A popular and simple example to accelerate the gradient method is by the Polak--Ribiere rule that calculates a `conjugate direction' as
\eq{eq:Pnew} 
    P = G  + \beta P_{\textrm{old}} \quad
    \text{with} \quad
    \beta = \frac{\left\langle G-G_{\textrm{old}}, G \right\rangle }
    {\left\langle G_{\textrm{old}}, G_{\textrm{old}}\right\rangle } . 
    \en
    Here, we avoid indices by calling $G_{\textrm{old}}$ the old gradient (usually indexed by $k$) and $G$ the new one
    (usually indexed by $k+1$).
    The inner product used above is the standard Frobenius inner product of matrices
    where $\left\langle X, Y \right\rangle = \Tr  ( Y^T X ) $. It is typical to restart with a pure gradient step ($\beta = 0$) when $P$ is not a descent direction and at every $k_{\rm restart}$ iterations for some fixed choice for $k_{\rm restart}$.
    
    When applied to objective functions defined on manifolds, two modifications are required to the Euclidean update in~\eqref{eq:Pnew}. First, since $G_{\textrm{old}}$ is a tangent vector of $X_{\textrm{old}}$, it needs to be `transported' to the current iterate $X$ in order for the inner product $\left\langle G_{\textrm{old}}, G \right\rangle$ to be well defined. A simple solution\footnote{It is known that this is a vector transport that is invariant to the choice of representative of the subspaces when the retraction on Grassmann is done via the polar factor, as we do in Alg.~\ref{alg:CGM}.} is by orthogonal projection onto the tangent space:
    \[
    \beta = \frac{\left\langle G-(I-XX^T) G_{\textrm{old}} , G \right\rangle }
    {\left\langle G_{\textrm{old}}, G_{\textrm{old}}\right\rangle } . 
    \]
    Since $G = (I-XX^T)G$, we do not need to compute this projection explicitly and the formula for $\beta$ in \eqref{eq:Pnew} remains valid in our case. Next, since $P$ is required to be a tangent vector, the result in \eqref{eq:Pnew} is again projected onto the tangent space as $(I - X X^T) P$.

\subsection{Line search} 
In order to use $P$ instead of $G$, we need to modify the line search in Algorithm~\ref{alg:gradA1}. We will explain the differences for a general $P$.

Let $X(\mu) = X_{k+1}$ and $X = X_k$ denote the new and old orthonormalized iterates. As before, we construct an iteration
\[
 X(\mu) = (X - \mu P) M
\]
where the search direction $P$ is a tangent vector, $P^T X = 0$, and gradient-related, $\Tr(G^T P) > 0$ with $G =\grad \phi(X)$.  In addition, $M$ is a normalization matrix such that $X(\mu)^T X(\mu) = I$.

A small calculation shows that the same normalization idea for $M$ from the gradient method (when $P=G$) can be used here: from the eigenvalue decomposition 
\[
V D_\beta V^T  = P^T P
\]
we define
\[ 
 D_\mu = (I + \mu^2 D_\beta)^{1/2}.
\]
Then it is easy to verify that
\begin{equation}\label{eq:Xnew_P}
 X(\mu) = (X - \mu P) V D_\mu^{-1} V^T
\end{equation}
has orthonormal columns and represents again the polar factor of $X - \mu P$.

Let $P_v = PV$ and $X_v = XV$. To perform the line search for $\mu$, we evaluate $\phi$ in the new point:
\begin{eqnarray}
  \phi (X(\mu))
  &=& -\tfrac{1}{2}  \Tr  (D_\mu\inv V^T (X-\mu P)^T A (X-\mu P) V D_\mu\inv ) \nonumber \\
  &=& -\tfrac{1}{2}  \Tr  (D_\mu\inv (X_v-\mu P_v)^T A (X_v-\mu P_v) D_\mu\inv ) \nonumber \\
  &=& -\tfrac{1}{2} \Tr  \left( D_\mu^{-2} \left( X_v^T A X_v - 2 \mu (P_v^T A X_v )  
     + \mu^2 ( P_v^T A P_v )\right) \right) \nonumber \\
  &=& -\tfrac{1}{2} \Tr   \left(   \left(I+\mu^2 D_{\beta} \right)^{-1} 
      \left( D_\alpha + 2 \mu \ D_\zeta 
        + \mu^2  \ D_\gamma \right) \right) \label{eq:new_phi_plus}
  \end{eqnarray}
  where
  \begin{equation}\label{eq:CG_linesearch_D}
  \begin{aligned}
   D_\alpha &= \Diag(X_v^T A X_v), &
   D_\beta &= \Diag(P_v^T P_v), \\ 
   D_\gamma &= \Diag(P_v^T A P_v), & 
   D_\zeta &= -\Diag(P_v^T A X_v).
  \end{aligned}
\end{equation}  
Comparing to \eqref{eq:ratmu1}, we see that a new $D_\zeta$ has appeared. Observe that $D_\alpha, D_\beta, D_\gamma$ all have non-negative diagonal but this is not guaranteed for $D_\zeta$. If $P=G$, then $-P_v^T A X_v = P_v^T P_v$ and thus $D_\zeta = D_\beta$. For a gradient related $P$ that is a tangent vector, we know that $0 \leq \Tr(P^T G) = -\Tr(VP^T \Pi AX V) =  -\Tr(P_v^T AX_v) = \Tr(D_\zeta)$. However, that does not mean that all the diagonal entries of $D_\zeta$ are non-negative, only their sum is. This lack of positive diagonal complicates the line search, as we will discuss next.

Let $\alpha_i, \beta_i, \gamma_i, \zeta_i$ be the $i$th diagonal entry of $D_\alpha, D_\beta, D_\gamma, D_\zeta$, resp. 
The rational function that represents~\eqref{eq:new_phi_plus} and generalizes~\eqref{eq:phinew} satisfies
\begin{equation}\label{eq:phi_P}
    \phi (X(\mu)) = - \half \sum_{i=1}^p  
    \frac{ \alpha_{i}\ + 2 \zeta_i  \mu +  \gamma_{i} \mu^2 }
    {1 +  \beta_i \mu^2} ,
\end{equation}
 with derivative
\begin{equation}\label{eq:deriv_P}
    \frac{ d \phi(X(\mu))}{d \mu} = -\sum_{i=1}^p \
    \frac{ \zeta_i\ + (\gamma_{i}\ -  \alpha_{i}\ \beta_i) \mu -      \beta_i \zeta_i\ \mu^2 }
    {(1+ \beta_i \mu^2)^2} \ .
\end{equation}
Since we do not know the sign of $\zeta_i$, each term in~\eqref{eq:deriv_P} has a quadratic in the numerator that can be convex or concave. This is different
from~\eqref{eq:dqnew}, where it is always convex (accounting for  the negative
sign outside the sum) since $\zeta_i=\beta_i$. In case there is a term with a concave quadratic, we can therefore not directly repeat the same arguments for the bracketing interval of $\mu$ based on the zeros of the quadratics in~\eqref{eq:deriv_P}. When there are negative $\zeta_i$'s, we could restart the iteration and replace $P$ by the gradient $G$. Since this wastes computational work, we prefer to simply disregard the branches that are concave when determining the bracket interval.

Overall, the line search for the CG approach will cost a little more than that for the gradient method, since we have an additional (diagonal) matrix to compute, namely $D_\zeta$.

\begin{algorithm}
  \caption{$\mbox{\texttt{Riemannian Conjugate Gradient Descent}}(A,X)$}\label{alg:CGM} 
  \begin{algorithmic}[1]
    \State  \textbf{Start:} Select initial $X_0$ such that $X_0^T X_0 = I$. Set $G=P=0$.
    
    \For{$k = 0, 1, \ldots$} 
    \State Keep $G_{\rm old} := G$.
    \State Update $G := \grad \phi(X_k) = -(A X_k - X_k C_k)$ with $C_k = X_k^T A X_k$.
    \If {$\| G \| < \textrm{tol}$}
    \State \Return
    \EndIf

    \State Diagonalize  $G^T G = V D_{\beta} V^T$.

    \State Compute $D_{\alpha}, D_{\gamma}, D_{\eta}$ from \nref{eq:CG_linesearch_D} with $X=X_k$.

    \State  Compute $
      \beta = \left\langle G -G_{\rm old}, G \right\rangle  / 
      \left\langle G_{\rm old}, G_{\rm old}\right\rangle $
  \State     Update $P := (I-X_k X_k^T) (G + \beta P)$ 
  \If {restart}
  \State $P:=G$
  \EndIf

    \State Compute $\mu$ as the minimizer of \nref{eq:phi_P} using a modified version \texttt{Get\_Mu}.

    \State Compute $X_{k+1}$ as the polar factor of $X_k - \mu P$ like in~\eqref{eq:Xnew_P}. 
\EndFor
\end{algorithmic}
\end{algorithm}


\section{Numerical implementation and experiments}

\subsection{Efficient and accurate implementation} \label{sec:numerical_implementation_details}

A proper numerical implementation of Algorithms~\ref{alg:gradA1} and~\ref{alg:CGM}, and in particular the line search, is critical to obtain highly accurate solutions. We highlight here four important aspects. 

In addition, we give some details on how to improve the efficiency of a direct implementation of these algorithms so that they require the same number of matrix vector products with $A$ as subspace iteration and LOBCG.

\paragraph{Calculation of bracket}
The $\beta_i$'s in~\eqref{eq:xii} can be very small in some situations. If we
set $\delta_i = \gamma_i - \alpha_i \beta_i $ then cancellation may cause
loss of accuracy in formula \eqref{eq:xii} 
when  $\delta_i <0$. We can circumvent this by observing that in this case:
\eq{eq:derv1} 
  \xi_{i}
  =   \frac{    \sqrt{  \delta_{i}^2 + 4 \beta_i^3 }  - |\delta_i|}{ 2 \beta_i^2 }
  =   \frac{  4 \beta_i^3 }{ 2 \beta_i^2 ( |\delta_i| + \sqrt{  \delta_{i}^2 + 4 \beta_i^3 } ) }
  =   \frac{  2 }{ | \delta_i/\beta_i| + \sqrt{  (\delta_{i}/\beta_i)^2 + 4 \beta_i } }. 
\en
 When $\delta_i>0$ we can simply use  \nref{eq:xii} which we rewrite as
\eq{eq:derv2}
  \xi_{i}
  =
 \frac{1}{2 \beta_i}
 \left(  \frac{\delta_i}{\beta_i}
 + \sqrt{  \left(\frac{\delta_{i}}{\beta_i}\right)^2 + 4 \beta_i } \right) . 
\en

\paragraph{Calculation of the minimizer}

For numerical reasons, it is advisable to compute a root of $\phi'$ instead of a minimum of $\phi$. This can be done in an effective way by a safe-guarded root finding algorithm, like the Dekker--Brent algorithm from \texttt{fzero} in \textsc{Matlab}. Since this algorithm converges superlinearly, we rarely need more than 10 function evaluations to calculate the minimizer of $\phi$ in double precision.

\paragraph{Efficient matvecs}

At each iteration $k$, the line search requires $AP_k$ and $AX_k$; see~\eqref{eq:CG_linesearch_D}. Supposing that $AX_k$ was calculated previously, it would seem that we need another multiplication of $A$ with $P_k$ which is not needed in subspace iteration (accelerated by Chebyshev or not). Fortunately, it is possible to avoid one of these multiplications. First, we proceed as usual by computing the next subspace $X_{k+1}$ from the polar decomposition
\[
 X_{\rm new} =  (X - \mu P) V D_\mu\inv V^T . 
\] 
Instead of calculating $A X_{\rm new}$ explicitly in the next iteration, we observe that
\begin{equation}\label{eq:save_matvec}
 A X_{\rm new} = (AX - \mu A P) V D_\mu\inv V^T.
\end{equation}
Hence, it suffices to compute only $A P$ explicitly at each iteration since $AX$ can be updated by the recursion above. 
Except for a small loss of accuracy when the method has nearly converged, this computation behaves very well numerically. In practice, the product $AX$ is only calculated explicitly when $\mu = O(\varepsilon_{\rm mach})$.

\paragraph{Efficient orthonormalization}

The line search procedure requires the diagonalization $P^T P = V D_\beta V^T$ which has a non-negligible cost of $O(np^2 + p^3)$ flops. Fortunately, the result of this decomposition can be used again for the normalization of $X_{new}$ by the polar factor, as explained in~\eqref{eq:Ynew} and~\eqref{eq:Xnew_P}. Compared to using QR for the normalization, there is therefore very little overhead involved.

\subsection{Comparison with subspace iteration for a Laplacian matrix}

We first test our methods for the standard 2D finite difference Laplacian on a $35 \times 40$ grid, resulting in a symmetric positive definite matrix of size $n = 1\,400$. Recall that the dimension of the dominant subspace to be computed is denoted by $p$.

Algorithms~\ref{alg:gradA1} and~\ref{alg:CGM} (with $k_{\rm restart} = 75$) are compared to subspace iteration applied to a shifted and scaled matrix $(A-c I) / h$ and a filtered matrix $p_d(A)$ with given degree $d$. The shift $c$ and scaling $h$ are defined in~\eqref{eq:optimal_poly}. Likewise, the polynomial $p_d$ is the one from~\eqref{eq:optimal_poly}, determined from a Chebyshev polynomial to filter the unwanted spectrum in $[\lambda_n, \ldots, \lambda_{p+1}]$. See also~\cite{chebfsi} for a concrete implementation based on a three-term recurrence that only requires computing one product $AX_k$ per iteration. Recall that these choices of the shift and the polynomial are in some sense optimal as explained \S\ref{sec:SI} for the given degree $d$. In addition, we compared to the locally optimal block conjugate gradients method (LOBCG) from~\cite{knyazev2001toward} which is closely related to Riemannian CG but with a higher cost per iteration; see \S\ref{sec:LOPBCG} for more details.

Observe that both subspace iteration methods make use of the exact values of the smallest eigenvalue $\lambda_n$ and of the largest unwanted eigenvalue $\lambda_{p+1}$. While this is not a realistic scenario in practice, the resulting convergence behavior should therefore be seen as the best case possible for those methods. Algorithms~\ref{alg:gradA1} and~\ref{alg:CGM} on the other hand, do not require any knowledge on the spectrum of $A$ and can be applied immediately.

The subspace iteration with Chebyshev acceleration will restart every $d$ iterations to perform a normalization of $X_k$ and, in practice, adjusts the Chebyshev polynomial based on refined Ritz values\footnote{This is not done in our numerical tests since we supply the method the exact unwanted spectrum.}. For small $d$, the method does not enjoy as much acceleration as for large $d$. On the other hand, for large $d$ the method is not stable. 

\begin{figure}
    \centering
    \includegraphics[width=0.75\textwidth]{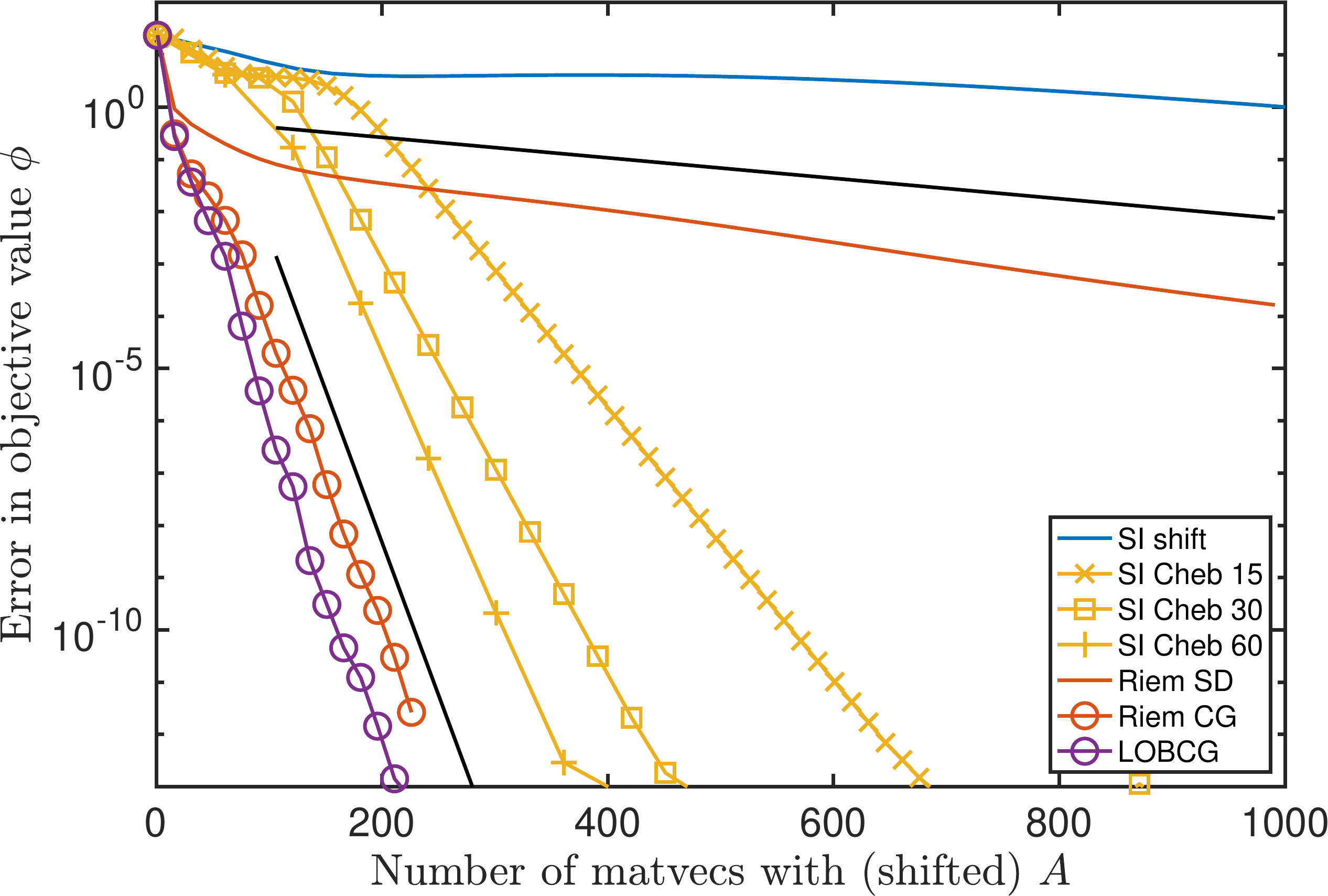}
    \caption{Error in objective value for subspace iteration (SI),  Riemannian steepest descent (SD), Riemannian nonlinear conjugate gradients (CG), and locally optimal block conjugate gradients (LOBCG) for a Laplacian matrix of size $n = 1\, 400$ based on finite differences when computing the dominant subspace of dimension $p=6$. For SI, optimal shift and optimal Chebyshev polynomials were used of various degree (number in legend).  The black lines estimate the asymptotic convergence speed as explained in the text.}
    \label{fig:Lapl35x40_trace}
\end{figure}

In Figure~\ref{fig:Lapl35x40_trace}, the convergence of the objective function $\phi(X_k)$ is visible for subspace dimension $p=6$ and polynomial degrees $d \in \{ 15, 30, 60\}$. All methods perform per iteration only one block matvec of $A$ with a matrix of size $n \times p$. Since this is the dominant cost in large-scale eigenvalue computations like SCF, we plotted the convergence in function of this number\footnote{For this example with very sparse $A$, the SI methods are much faster per iteration than the Riemannian methods. This is mainly because SI only needs to orthonomalize every $d$ times.}.

The benefits of acceleration by the Chebyshev polynomial filter or by Riemannian CG are clearly visible in the figure. In black lines, we also indicated the asymptotic convergence $O(\gamma^k)$ in function of the number of matvecs $k$ for two values of $\gamma$. In particular, it is well known (see, e.g., \cite[Lemma 17]{alimisis2022geodesic}) that
\begin{equation}\label{eq:kappa}
 \kappa = \frac{\lambda_1 - \lambda_n}{\lambda_p - \lambda_{p+1}} = \mathcal{O}(1/\delta).
\end{equation}
is the condition number of the Riemannian Hessian of $\phi$ at the dominant subspace $\mathcal{V}_\alpha$ with spectral gap $\delta$. From this, the asymptotic convergence rate of Riemannian SD is known (see \cite[Chap.~12.5]{luenbergerLinearNonlinearProgramming2008}) to satisfy
\[
 \gamma_{SD} = \left( \frac{\kappa - 1}{\kappa + 1} \right)^2 = 1 - \mathcal{O}(\delta).
\]
In addition, for Riemannian CG we conjecture the rate
\[
 \gamma_{CG} = \left( \frac{\sqrt{\kappa} - 1}{\sqrt{\kappa} + 1} \right)^2 = 1 - \mathcal{O}(\sqrt{\delta})
\]
based on the similarity to classical CG for a quadratic objective function with condition number $\kappa$. For both Algorithms~\ref{alg:gradA1} and~\ref{alg:CGM}, we see that the actual convergence is very well predicted by the estimates above. 

\subsection{A few other matrices}\label{sec:other_matrices}

As our next experiment, we apply the same algorithms from the previous section (but without restarting to have parameter free Riemannian methods) to a few different matrices
 and several choices for the subspace dimension $p$. In addition, we target also the minimal eigenvalues by applying the methods to $-A$ instead of $A$. This is not a problem, as the Riemannian gradient of $\phi$ on Grassmann is invariant under shifts. More concretely, Riemannian steepest descent (RSD) and Riemannian CG (RCG) with exact line search applied to $-A$ produce the same iterates as when applied to $-A+c I$, for any $c \in \mathbb{R}$. For Algorithm \ref{alg:CGM} the signs of $G$ and $G_{old}$ flip, but the parameter $\beta$ remains the same at each iteration. Thus, both methods converge to the eigenvectors associated to the largest eigenvalues of $-A$, which are the eigenvectors associated with the smallest eigenvalues of $A$.
 
 Except for the standard finite difference matrices for the 3D Laplacian, the matrices used were taken from the SuiteSparse Matrix Collection~\cite{SuiteSparse}. This results in problems with moderately large Riemannian condition numbers $\kappa$, defined in~\eqref{eq:kappa}.

Due to the larger size of some of these matrices, we first compute with a Krylov--Schur method (implemented in \textsc{Matlab} as \texttt{eigs}) the eigenvalues that are required to determine the optimal Chebyshev filter in  subspace iteration. The Riemannian methods do not require this or any other information. As optimal value $\phi^*$ for the function value, we took the best value of the results computed from all methods, including the Krylov--Schur method.  

\paragraph{FD3D} This matrix is the 3D analogue of the matrix we tested in the previous section. It corresponds to a standard finite difference discretization of the Laplacian in a box with zero Dirichlet boundary conditions. We used $n_x=35, n_y = 40, n_z=25$ points in the $x,y,z$ direction, resp. The resulting matrix is of size $35\, 000$. Compared to the earlier experiment, we took larger subspace dimensions and also a minimization of the Rayleigh quotient. All these elements make for a more challenging problem numerically.

\medskip

\begin{center}
\begin{tabular}{ c c c c c } \toprule
 problem  & type & dimension  $p$   & Riem.~cond.~nb. & Cheb.~degree \\\midrule
 1 & min & 64   & $3.53 \cdot 10^4$ & 100 \\
 2 & max & 32  & $5.54 \cdot 10^3$ & 100 \\ \bottomrule
 \end{tabular}
\end{center}

\medskip 

\begin{figure}
\begin{minipage}{0.48\textwidth}
  \includegraphics[width=0.99\textwidth]{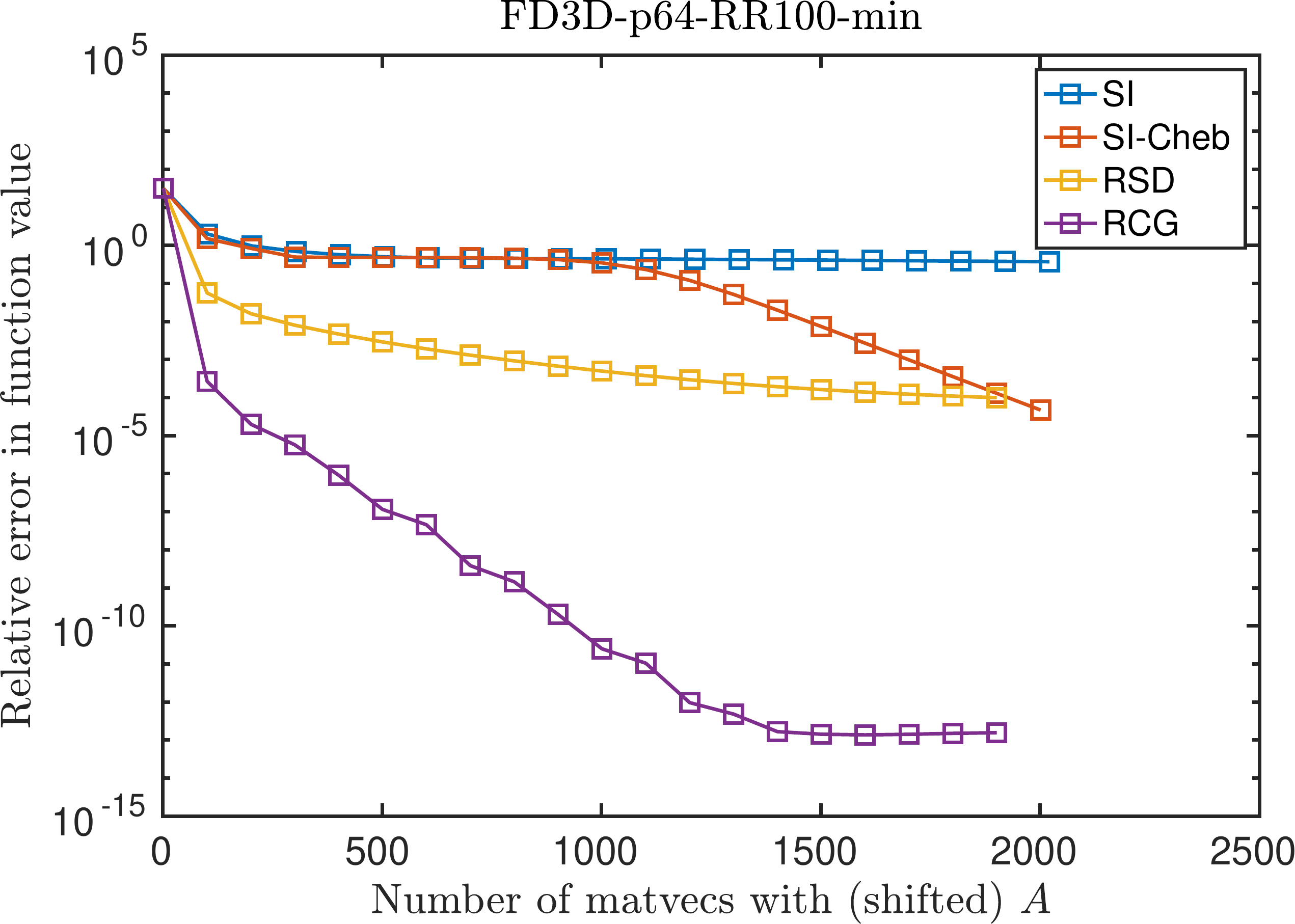} 
   \centerline{Problem nb.~1}
  \end{minipage}
  \begin{minipage}{0.48\textwidth}
  \includegraphics[width=0.99\textwidth]{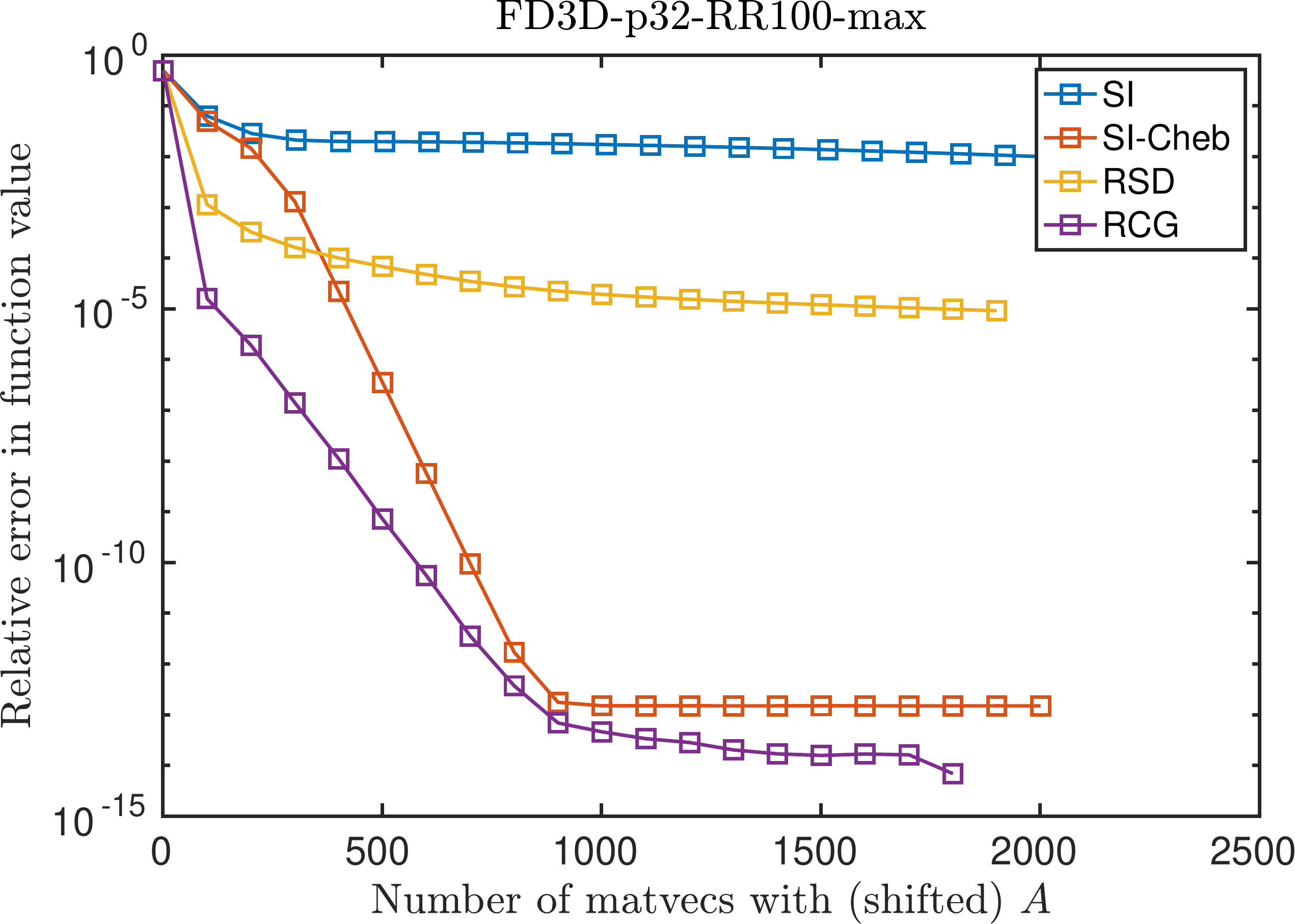} 
    \centerline{Problem nb.~2}
\end{minipage} 
\caption{The FD3D matrix.}\label{fig:fd3d}   
\end{figure}

In Fig.~\ref{fig:fd3d}, we see that the convergence of the maximization problem is very similar to that of the 2D case, although the asymptotic convergence rate of Riemannian CG seems to be slower than that of subspace iteration with optimal filter. This can be improved by restarting (not shown) but even without it, the results are good. On the other hand, the more relevant case of finding the minimal eigenvalues of a Laplacian matrix turns out to be a challenge for SI with or without Chebyshev acceleration. In fact, even with a degree 100 polynomial it takes about 1000 iterations before we see any acceleration. The Riemannian methods, on the other hand, converge much faster and already from the first iterations. 

\bigskip 

\paragraph{ukerbe1} This matrix
is related to a 2D finite element problem on a locally refined grid
and it has a  relatively small size $n=5\,981$.
It is therefore more interesting than the uniform grid of the Laplacian
examples above. We tested the following parameters.

\medskip 

\begin{center}
\begin{tabular}{ c c c c c } \toprule
 problem  & type & dimension  $p$   & Riem.~cond.~nb. & Cheb.~degree \\\midrule
 3 & max & 32   & $4.85 \cdot 10^3$ & 50 \\
 4 & max & 64   & $5.21 \cdot 10^3$ & 100 \\ \bottomrule
 \end{tabular}
\end{center}

\medskip 

\begin{figure}
\begin{minipage}{0.48\textwidth}
  \includegraphics[width=0.99\textwidth]{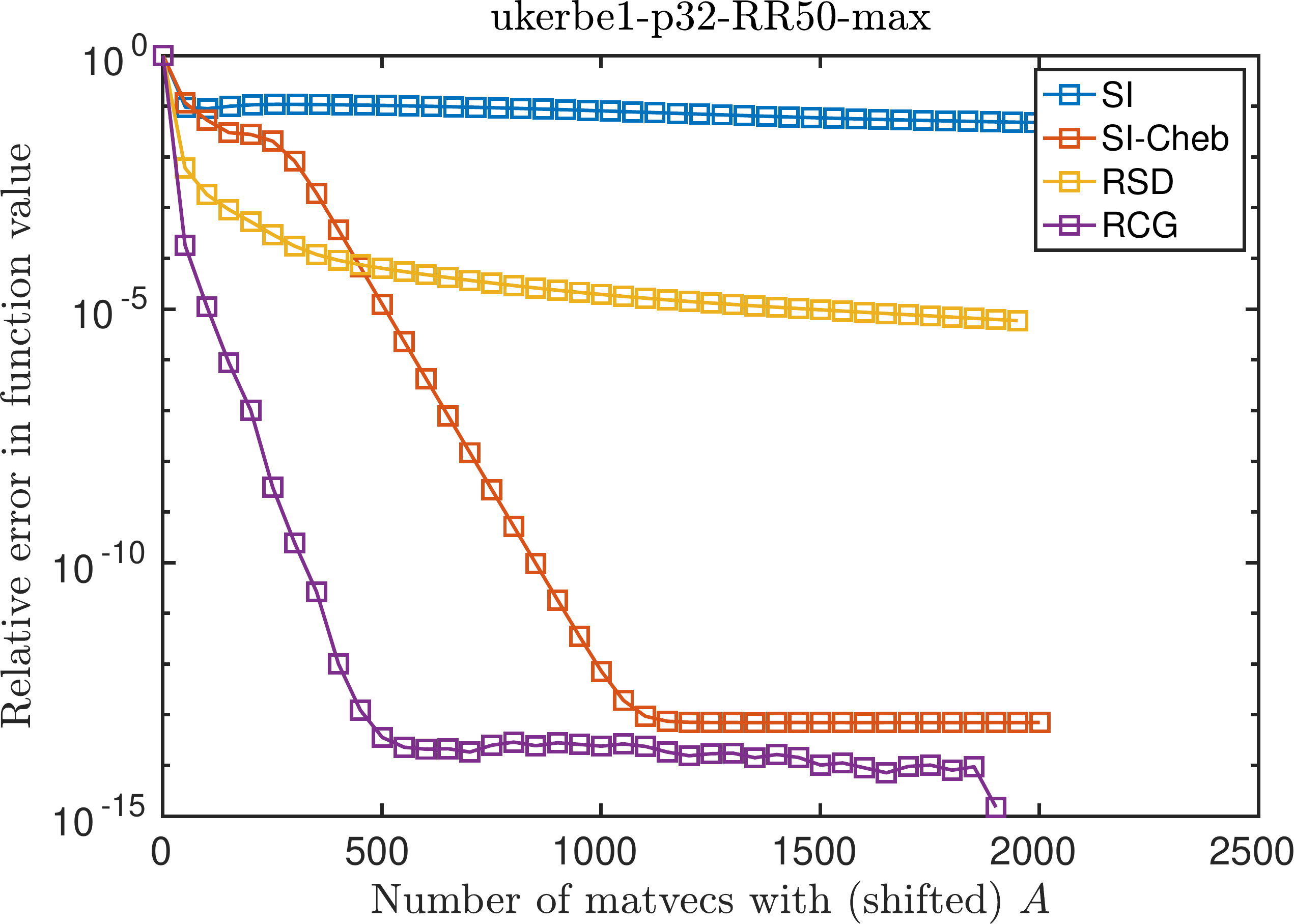} 
   \centerline{Problem nb.~3}
  \end{minipage}
  \begin{minipage}{0.48\textwidth}
  \includegraphics[width=0.99\textwidth]{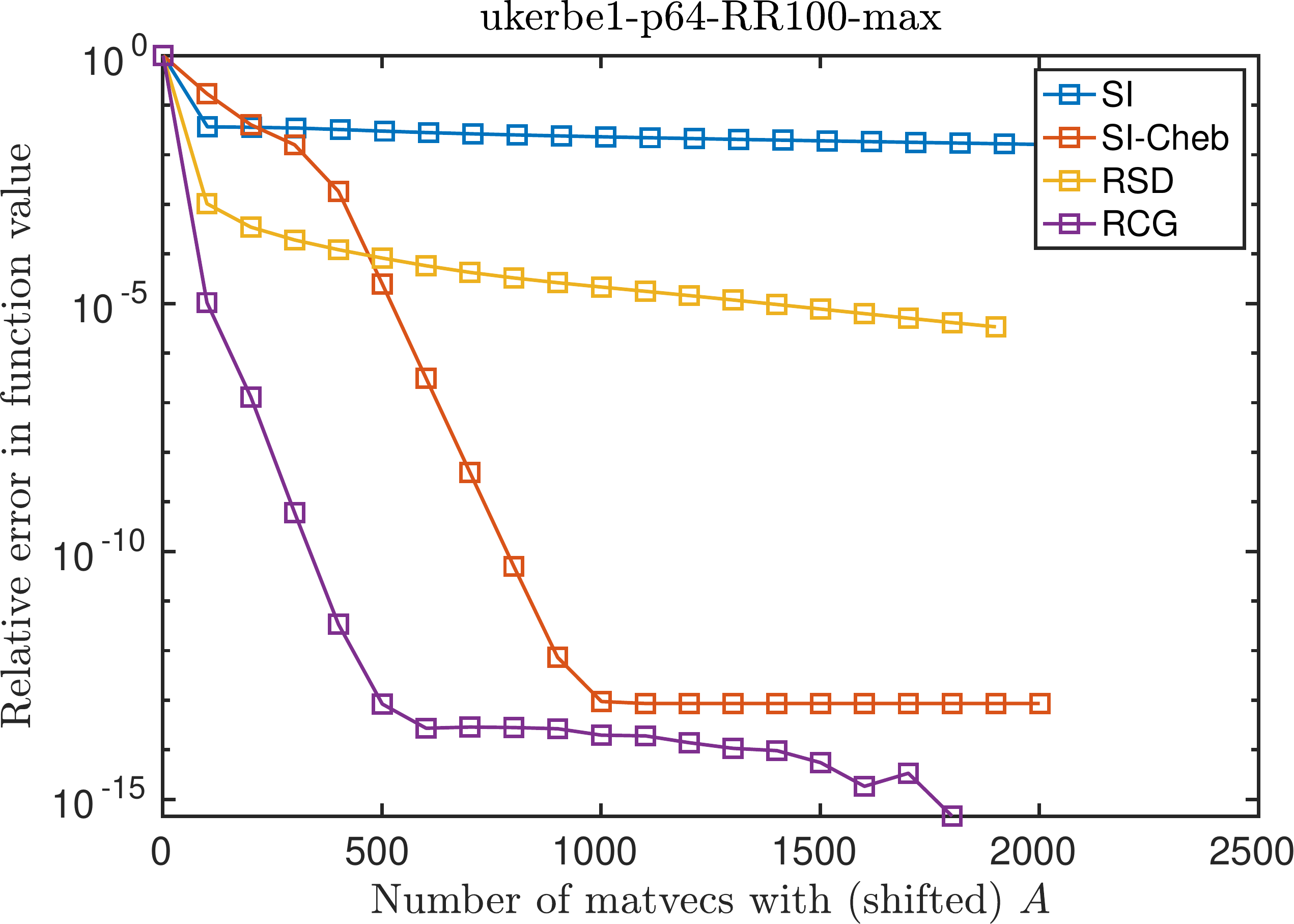} 
    \centerline{Problem nb.~4}
\end{minipage} 
\caption{The ukerbe1 matrix.}\label{fig:ukerbe}   
\end{figure}

In Figure~\ref{fig:ukerbe}, we observe that the Riemannian algorithms converge faster than their subspace iteration counterparts. This behavior is seen for many choices of $p$ and the Chebyshev degree. Since the spectrum of this matrix is symmetric around zero, the min problems are mathematically equivalent to the max problems, and therefore omitted.

\bigskip

\paragraph{ACTIVSg70K} We now test a larger matrix of size 69\,999. It models a synthetic (yet realistic) power system grid from the Texas A\&M Smart Grid Center. This matrix has a spectral gap of $\mathcal{O}(10)$ but the Riemannian condition number, which represents the correct relative measure of difficulty, is still large. Such a different kind of scale makes this an interesting matrix to test our algorithms.

\medskip 

\begin{center}
\begin{tabular}{ c c c c c } \toprule
 problem  & type & dimension  $p$  & Riem.~cond.~nb. & Cheb.~degree \\\midrule
 5 & min & 16 & $1.15 \cdot 10^4$ & 50 \\
 6 & max & 32 & $1.29 \cdot 10^3$ & 50 \\ \bottomrule
 \end{tabular}
\end{center}

\medskip 

For the minimization problem (nb.~5), we see that the Riemannian algorithms converge considerably faster than subspace iteration with or without Chebyshev acceleration of degree 50. (The reason for the bad performance of the Chebyshev acceleration is due to numerical instability with a degree 50 polynomial for this problem.) 
For the maximization problem (nb.~6), Riemannian CG and Chebyshev acceleration with degree 50 have very similar asymptotic convergence speed although the Riemannian algorithm has a faster start. The same conclusion hods for Riemannian SD and standard subspace iteration, although their convergence is of course significantly slower.

\begin{figure}
\begin{minipage}{0.48\textwidth}
  \includegraphics[width=0.99\textwidth]{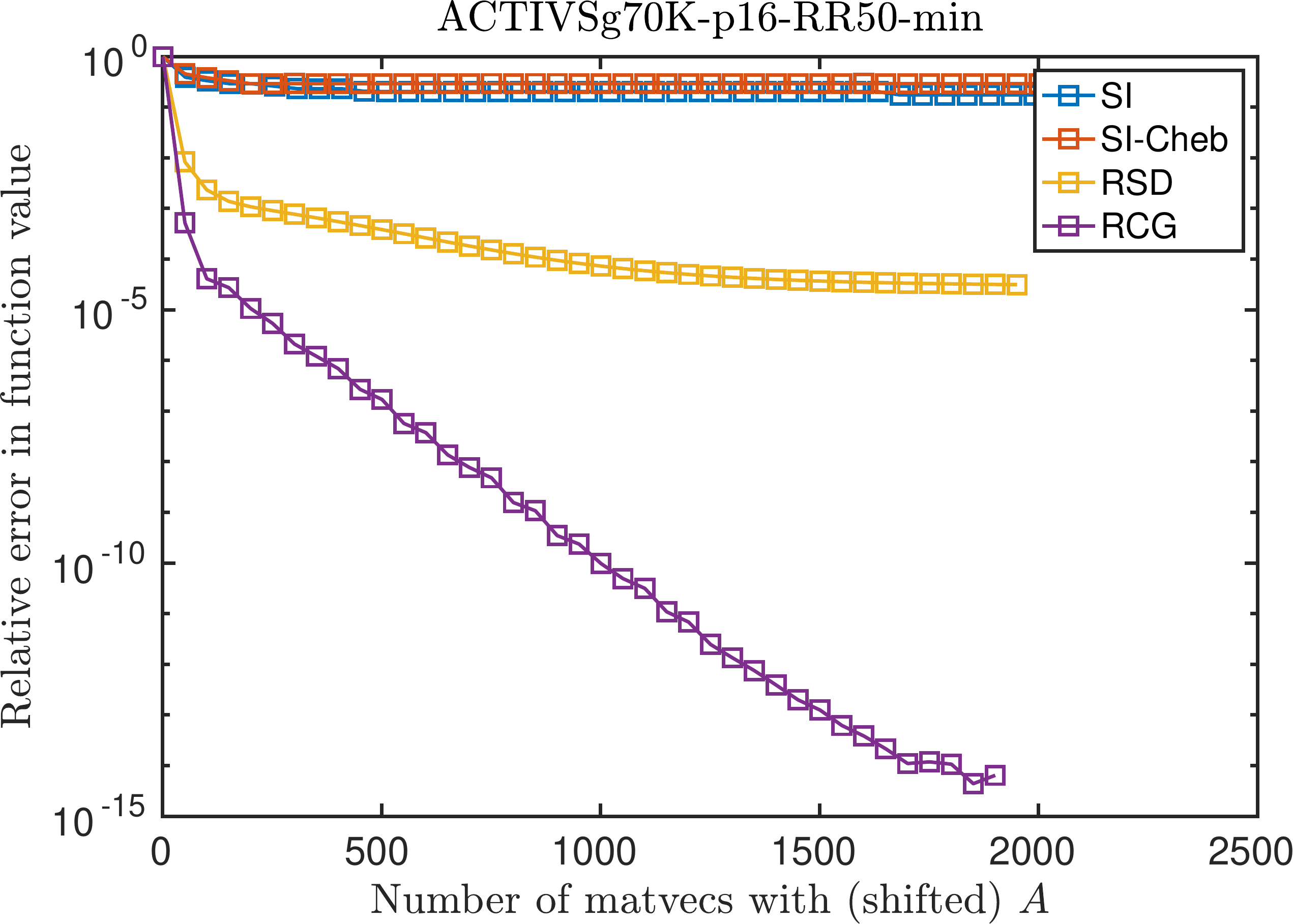} 
    \centerline{Problem nb.~5}
  \end{minipage}
  \begin{minipage}{0.48\textwidth}
  \includegraphics[width=0.99\textwidth]{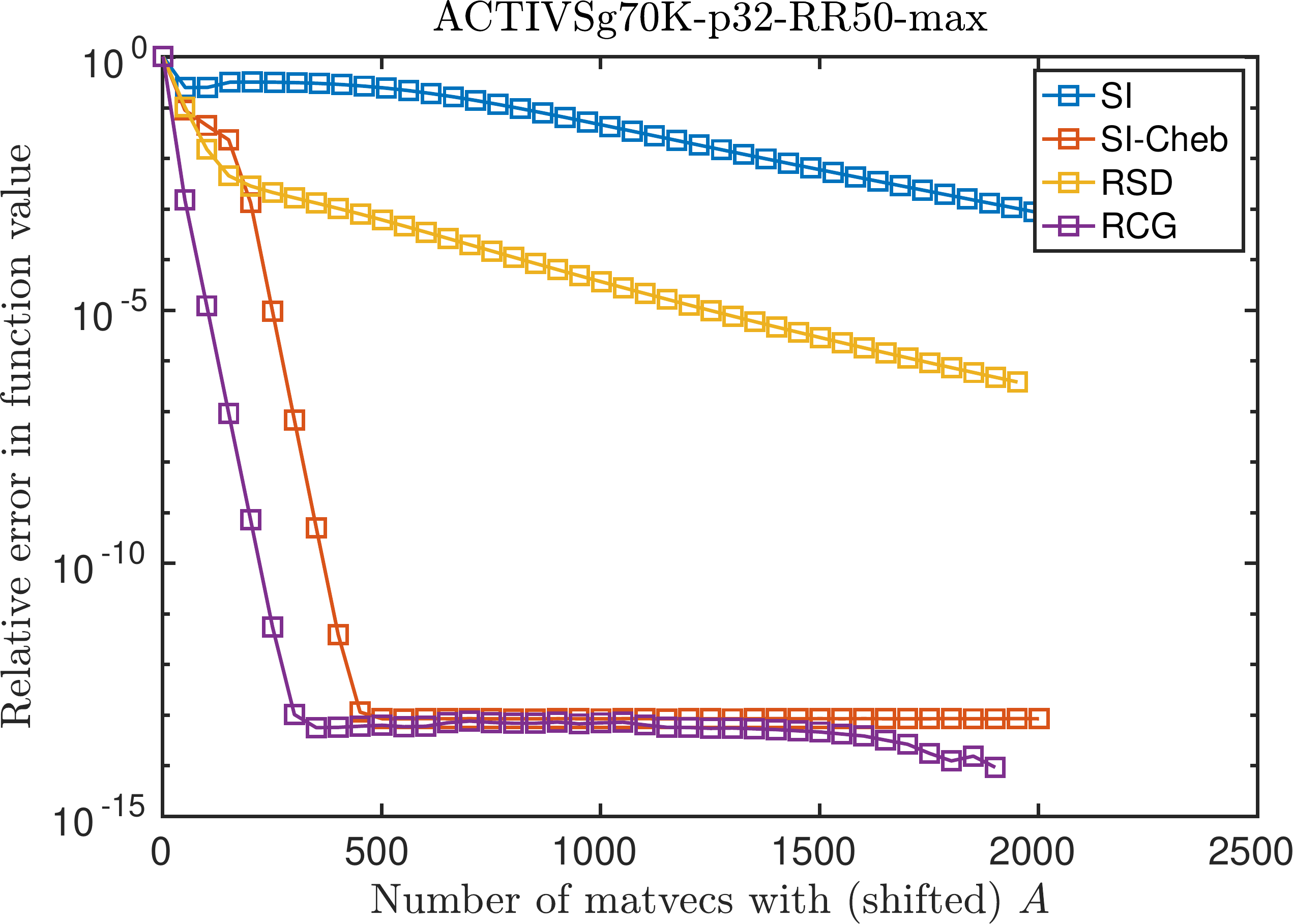} 
    \centerline{Problem nb.~6}
\end{minipage} 
\caption{The ACTIVSg70K matrix.}\label{fig:ACTIVS}   
\end{figure}

\bigskip 

\paragraph{boneS01} This final matrix is part of the Oberwolfach model order reduction benchmark set and models a 3D trabecular bone. It is our largest example of size $n=127\,224$. As we can see from the table below, for subspace dimension $p=64$ the minimization problem is particurlay challenging with a large Riemannian condition number.

\medskip 

\begin{center}
\begin{tabular}{ c c c c c } \toprule
 problem  & type & dimension  $p$   & Riem.~cond.~nb. & Cheb.~degree \\\midrule
 7 & min & 64  & $2.57 \cdot 10^6$ & 25 \\
 8 & max & 64  & $2.05 \cdot 10^3$ & 25 \\ \bottomrule
 \end{tabular}
\end{center}

\medskip 

\begin{figure}
\begin{minipage}{0.48\textwidth}
  \includegraphics[width=0.99\textwidth]{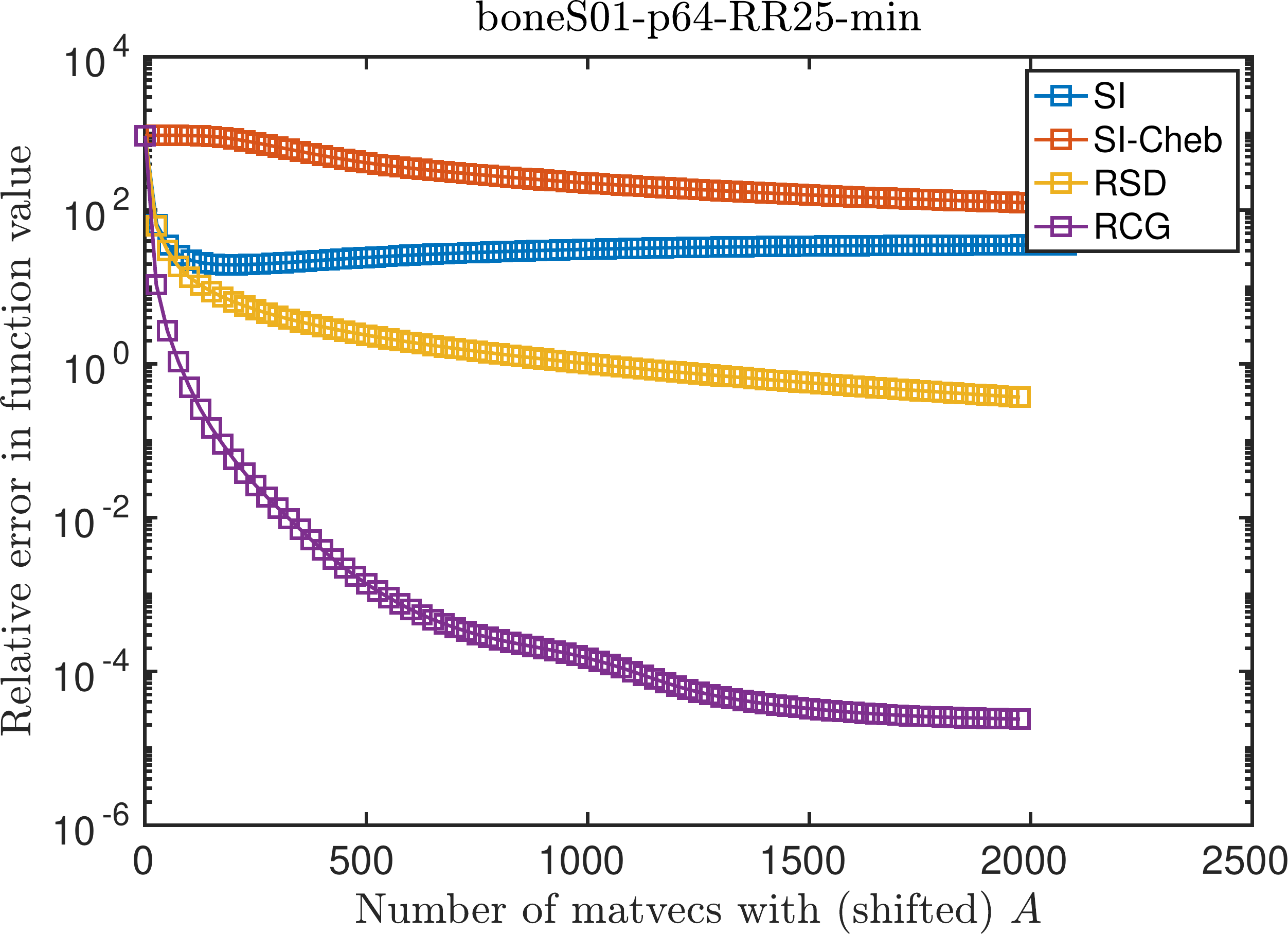} 
    \centerline{Problem nb.~7}
  \end{minipage}
  \begin{minipage}{0.48\textwidth}
  \includegraphics[width=0.99\textwidth]{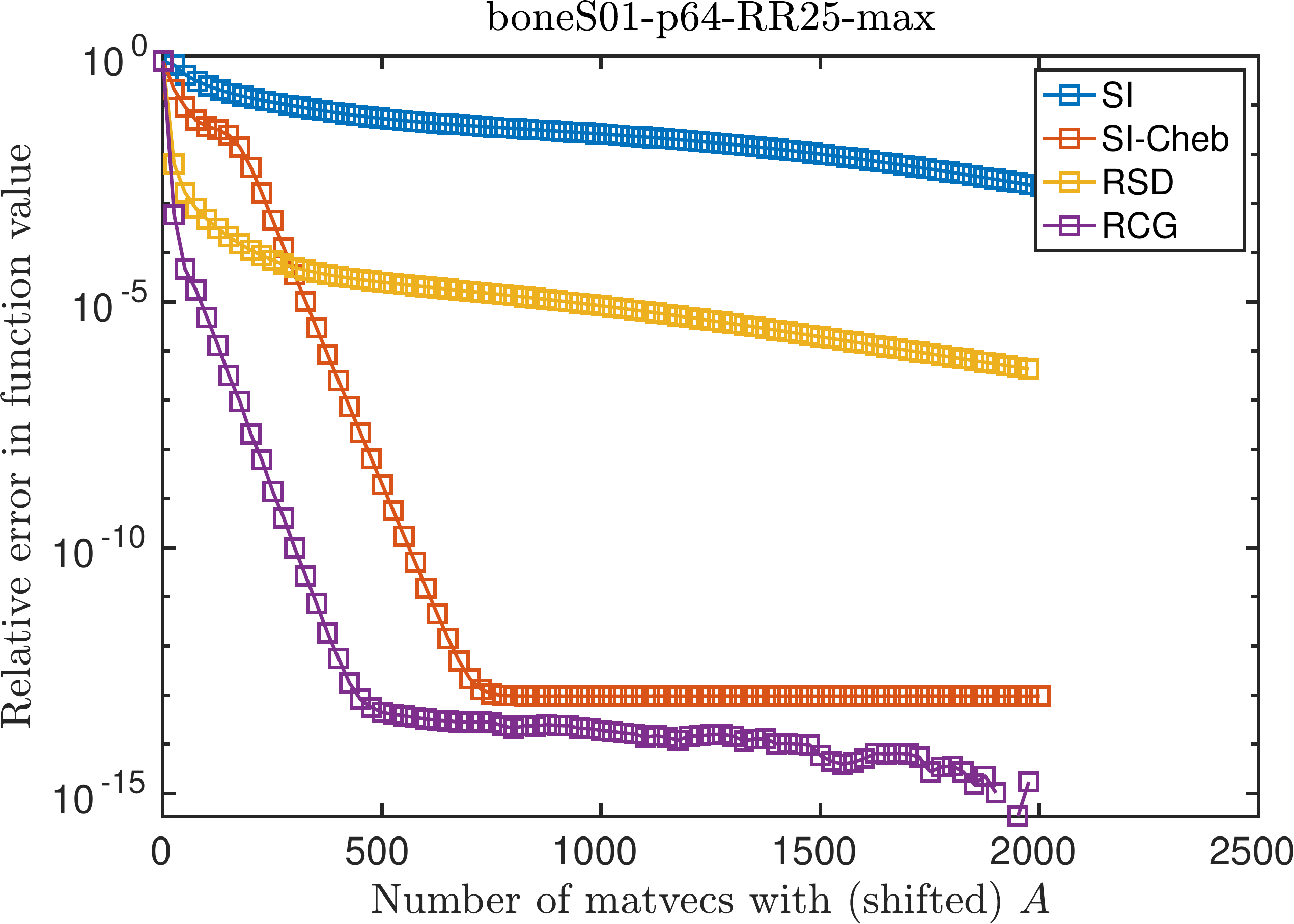} 
    \centerline{Problem nb.~8}
\end{minipage} 
\caption{The boneS01 matrix.}\label{fig:bone}   
\end{figure}

The convergence of the methods is visible in Fig.~\ref{fig:bone}. We can make similar observations as for the example above: the Riemannian algorithms have a faster initial convergence compared to the subspace variants. In addition, the accelerated variants are clear improvements.

\bigskip

\subsection{Comparison to LOBCG} \label{sec:LOPBCG}

It is instructive to compare the Riemannian CG method from Alg.~\ref{alg:CGM} to the locally optimal block CG method (LOBCG) from~\cite{knyazev2001toward} since both methods minimize the partial trace function $\phi$ using momentum terms. LOBCG is equivalent to the better known LOPBCG method where the preconditioner is not used.

Let $k$ be the iteration number. The essential difference between the two methods is that LOBCG minimizes $\phi$ over all orthonormal matrices that lie in the $3p^2$-dimensional subspace\footnote{When $X_k$ converges, adding $X_{k-1}$ to the columns of $X_k$ and $G_k$ would lead to numerical cancellation when computing an orthonormal basis for $\mathcal{V}_k$. In the implementation of LOBCG, a different matrix is therefore added that has better numerical properties. For theoretical investigations, we can ignore it.}  
\begin{equation}\label{eq:def_Vk_LOBPCG}
\mathcal{V}_k = \myspan(X_k, G_k, X_{k-1}) = \{ X_k \Omega + G_k \Psi + X_{k-1} \Phi \colon \Omega, \Psi, \Phi \in \mathbb{R}^{p \times p} \}.
\end{equation}
Here, the residual $G_k = AX_k - X_k X_k^T A X_k$ is also the Riemannian gradient of $\phi$ at $X_k$. Contrary to most optimization problems, this subspace search can be computed exactly for the symmetric eigenvalue problem by the Rayleigh--Ritz procedure: the optimal solution is related to the top $p$ eigenvectors of the symmetric $3p \times 3p$ matrix $Q_k^T A Q_k$ with $Q_k$ an orthonormal basis for $\mathcal{V}_k$. 

In contrast, the Riemannian CG method minimizes $\phi$ for the scalar $\alpha$  during the line search applied to the orthonormalization of $X_k - \alpha P_k$. When $p>1$, there is no explicit solution for the optimal $\alpha$ in terms of a smaller eigenvalue problem, but as explained above, it can be solved efficiently by diagonalizing the matrix $P_k^T P_k$.

When started at the same $X_k$ and $X_{k-1}$, LOBCG will produce $X_{k+1}$ with a smaller objective value $\phi(X_{k+1})$ than the Riemannian CG method. This is because an iterate produced with the step $X_k - \alpha P_k$ from Riemannian CG is contained in the subspace searched by LOBCG; see~App.~\ref{app:A} for a proof. It is therefore reasonable to expect\footnote{Since the iteration is not stationary and depends on the previous iterates, one cannot conclude that LOBCG always produces iterates with lower objective value than Riemannian CG.} that LOBCG converges faster overall in terms of number of iterations. 

In Table~\ref{tab:LOPBCG}, we have compared LOBCG to Riemannian CG (denoted by RCG) for the same matrices we tested above. For the matrices ukerbe1 and FD3D, we see that LOBCG indeed requires less iterations than Riemannian CG, usually by about a factor two. However, this does not mean that LOBCG is faster in terms of computational time due to an increased cost per iteration. In addition, the differences between LOBCG and Riemannian CG are less predictable for the other matrices. Overall, Riemannian CG is usually faster in computational time and also more reliable.

The increased cost per iteration of LOBCG compared to Riemannian CG is due to the additional computations for the subspace search. While both methods only require one product of the form $AZ$ with an $n \times p$ matrix $Z$, LOBCG performs 3 orthonormalizations (by Cholesky decomposition) whereas Riemannian CG needs 2 (by polar factor). Furthermore, LOBCG needs 14 matrix products of the form $Y^T Z$ for $n \times p$ matrices $Y$ and $Z$, while Riemannian CG requires only 4. Finally, the calculation of $X_{k+1}$ (and $AX_{k+1}$) based on the coefficients from the Rayleigh--Ritz procedure is not negligible in LOBCG with a cost comparable to a product $Y^T Z$. For Riemannian CG, it is simply a linear combination of two matrices (before normalization). In our experiments, one iteration of LOPBCG was therefore about 2 to 3 times more expensive, depending on $A$ and $p$.

\begin{figure}
\begin{minipage}{0.48\textwidth}
  \includegraphics[width=0.99\textwidth]{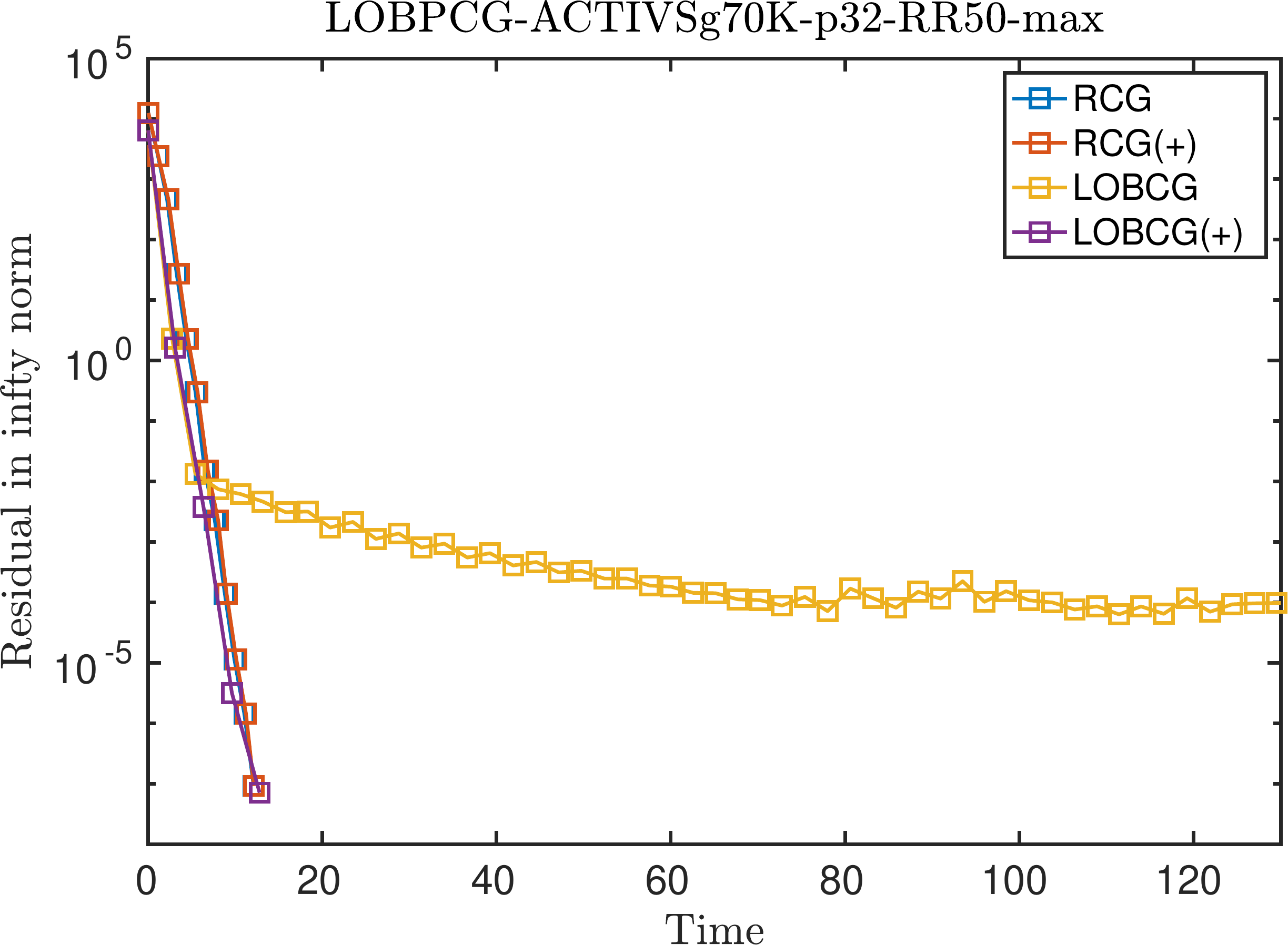}
  \end{minipage}
  \begin{minipage}{0.48\textwidth}
  \includegraphics[width=0.99\textwidth]{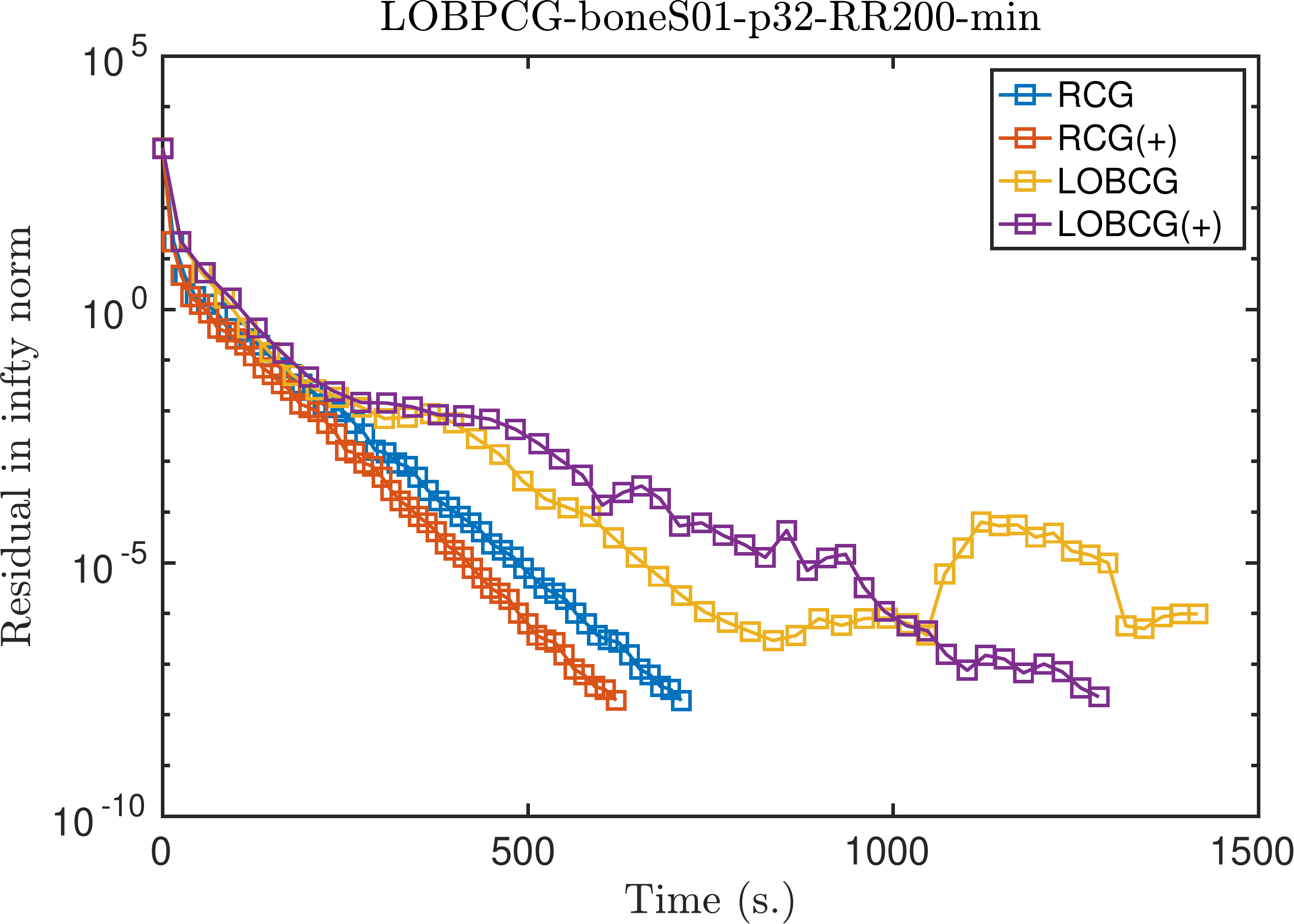} 
\end{minipage} 
\caption{Instability of the original LOBCG method.}\label{fig:unstable_LOBCG}   
\end{figure}

We have also tested a version of LOBCG where all the block entries in $Q_k^T A Q_k$ are explicitly calculated (denoted by LOBCG(+) in the table). The original code replaces $X_k^T A X_k$ by the eigenvalues obtained in the Rayleigh--Ritz procedure. While this behaves well early on, we have noticed stability issues in our experiments. Figure~\ref{fig:unstable_LOBCG}  is a clear example where the original version LOBCG does not converge or behaves erratically. In other examples (not shown), the residual even grows in an unbounded way. The version LOBCG(+) is however not always an improvement over LOBCG, which can be seen from the table. This shows that an accurate implementation of CG-based methods is not trivial, even with subspace search.

For Riemannian CG, we also tested a version (denoted by RCG(+)) where the product of $AX_k$ is explicitly calculated instead of being computed recursively as in~\eqref{eq:save_matvec}. The unchanged number of iterations in Table~\ref{tab:LOPBCG} shows that there is no loss of accuracy when utilizing the recursion. When the matrix $A$ is very sparse, like in FD3D, the version RCG(+) is less costly per iteration but for other matrices, the original version RCG is preferable.

\begin{table} \label{tab:LOPBCG}
\caption{Comparison of LOBCG and Riemannian CG (denoted by RCG) when minimizing/maximizing the partial trace for a few test matrices with different block sizes $p$. The time in seconds (rounded to nearest integer) and number of iterations to reach a relative residual $\|G_k\|_\infty / \|G_0\|_\infty$ of $10^{-8}$ is indicated in sec.~and its., resp. If the method did not reach the required tolerance in 10\,000 iterations, a star * is given. The methods indicated with a (+) are variants that aim to be more accurate; see the text for their definition.}
\begin{center}
\begin{tabular}{c c c c c c @{\hskip 0.8cm} c c c c} \toprule
 & & \multicolumn{2}{c}{LOBCG} & \multicolumn{2}{c@{\hskip 0.8cm}}{LOBCG(+)} &  \multicolumn{2}{c}{RCG(+)} & \multicolumn{2}{c}{RCG} \\
& problem  & secs. & its. & secs. & its. & secs. & its. & secs. & its. \\\midrule
\parbox[t]{2mm}{\multirow{6}{*}{\rotatebox[origin=b]{90}{ACTIVSg70K}}}     & p=16 max &    *   &   *    &     5 &  152  &     4 &  301  &     $\bm{3}$ &  301  \\
& p=16 min &     $\bm{4}$ &  152  &     $\bm{4}$ &  152  &    29 & 2151  &    25 & 2151  \\
 & p=32 max &   110 & 2152  &    10 &  152  &     $\bm{9}$ &  451  &    11 &  451  \\
& p=32 min &     5 &  102  &    *   &   *    &     $\bm{4}$ &  201  &     5 &  201  \\
 & p=64 max &    *   &   *    &    *   &   *    &    $\bm{70}$ & 1301  &    86 & 1301  \\
 & p=64 min &   401 & 2852  &    $\bm{19}$ &  102  &    27 &  551  &    36 &  551  \\ \midrule
 & p=16 max &     $\bm{9}$ &  752  &    10 &  752  &    11 & 1801  &    10 & 1801  \\
 & p=16 min &     7 &  602  &     8 &  652  &     3 &  451  &     $\bm{2}$ &  451  \\
\parbox[t]{2mm}{\multirow{3}{*}{\rotatebox[origin=c]{90}{FD3D}}}  & p=32 max &    16 &  552  &    16 &  552  &    $\bm{11}$ & 1051  &    $\bm{11}$ & 1051  \\
& p=32 min &    22 &  802  &   $\bm{20}$ &  702  &    37 & 3701  &    40 & 3701  \\
 & p=64 max &    51 &  752  &    57 &  802  &    $\bm{22}$ &  901  &    27 &  901  \\
 & p=64 min &    55 &  802  &    52 &  752  &    $\bm{36}$ & 1401  &    42 & 1401  \\ \midrule
 & p=16 max &    $\bm{25}$ &  352  &    27 &  352  &    39 &  501  &    $\bm{25}$ &  501  \\
 & p=16 min &    276 & 4202  &   287 & 4002  &   340 & 4401  &   $\bm{209}$ & 4401   \\
\parbox[t]{2mm}{\multirow{3}{*}{\rotatebox[origin=c]{90}{boneS01}}} & p=32 max &    42 &  252  &    52 &  302  &    22 &  301  &    $\bm{19}$ &  301  \\
 & p=32 min &    648 & 4202  &   825 & 5202  &   480 & 6601  &   $\bm{412}$ & 6601    \\
 & p=64 max &    *   &   *    &   170 &  402  &   $\bm{101}$ &  651  &   $\bm{101}$ &  651  \\
 & p=64 min &    *   &   *    &    *   &   *    &    *   &   *    &    *   &   *    \\  \midrule
& p=16 max &     $\bm{1}$ &  452  &     $\bm{1}$ &  502  &     $\bm{1}$ &  601  &     $\bm{1}$ &  601  \\
& p=16 min &     $\bm{1}$ &  452  &     $\bm{1}$ &  452  &     $\bm{1}$ &  501  &     $\bm{1}$ &  501  \\
\parbox[t]{2mm}{\multirow{3}{*}{\rotatebox[origin=c]{90}{ukerbe1}}}  & p=32 max &     2 &  552  &     2 &  502  &     $\bm{1}$ &  651  &     2 &  651  \\
& p=32 min &     2 &  402  &     2 &  402  &     2 &  701  &     $\bm{1}$ &  701  \\
& p=64 max &     4 &  352  &     4 &  352  &     $\bm{3}$ &  651  &     4 &  651  \\
& p=64 min &     5 &  502  &     5 &  452  &     $\bm{3}$ &  551  &     $\bm{3}$ &  551  \\ \bottomrule
 \end{tabular}
\end{center}
\end{table}

\section{Conclusion}

We revisited the standard Riemannian gradient descent method for the symmetric eigenvalue problem as a more competitive alternative of subspace iteration. If accelerated using a momentum term from nonlinear CG, there is a wide variety of matrices where the Riemannian method is faster per number of matrix vector products than subspace iteration with optimal Chebyshev filter polynomials and faster in computational time than LOBCG. This property would make it valuable in applications like the self-consistent field (SCF) iteration.

Among novel contributions, we derived a computationally efficient exact line search. Its accurate implementation is key to the good performance of the method. We also presented new convergence proofs for this geodesic-free Riemannian algorithm,
including  a locally fast convergence result  in a $\mathcal{O}(\sqrt{\delta})$ neighbourhood of the dominant subspace.

\section*{Acknowledgments}

FA was supported by SNSF grant 192363. YS was supported by NSF grant DMS-2011324. BV was supported SNSF grant 192129.

\bibliographystyle{siam}

\bibliography{strings,eig,saad,local,refs}

\appendix

\section{Suboptimality of Riemannian CG in the LOBCG subspace}\label{app:A}

We prove that Riemannian CG with $p\geq 1$ is suboptimal compared to LOBCG when started at the same $X_k$ and $X_{k-1}$. The case  $p=1$ is also explained in~\cite[Sections~4.6.5 and~8.3]{absilOptBook08}. This improvement is of course more computationally expensive.

Since Riemannian CG produces iterates of the form
\begin{equation}\label{eq:iterate_RCG}
 X_{k+1} = (X_k - \alpha_k P_k) M_k
\end{equation}
with $M_k$ the normalization so that $X_{k+1}$ has orthonormal columns, it is clear that
\[
 X_{k+1} \in \myspan(X_k, P_k).
\]
Here, $\myspan(\cdot,\cdot)$ is to be interpreted as in~\eqref{eq:def_Vk_LOBPCG} as a subspace of dimension $2p^2$. 
Since $P_k = (I-X_k X_k^ T) (G_k + \beta_k P_{k-1})$, we also have $P_k \in \myspan(G_{k},P_{k-1}, X_k)$, from which
\[
 X_{k+1} \in \myspan(X_k, G_k, P_{k-1}).
\]
The relation~\eqref{eq:iterate_RCG} also shows that $P_{k-1} \in \myspan(X_{k-1},X_k)$ if $M_{k-1}$ is invertible, which is true generically. We therefore get that 
\[
 X_{k+1} \in \myspan(X_k, G_k, X_{k-1}) = \mathcal{V}_k,
\] 
where $\mathcal{V}_k$ is the subspace used in LOBCG. Since LOBCG is optimal for $\phi$ over all orthonormal matrices with $p$ columns in $\mathcal{V}_k$, it will be a lower bound of $\phi(X_{k+1})$.

\end{document}